\documentclass[journal,onecolumn,12pt]{IEEEtran}
\usepackage{amssymb}
\usepackage{graphicx,times,psfig,amsmath,algorithmicx,algpseudocode,booktabs}
\usepackage{amsmath,amsfonts,bm,mathtools}
\usepackage{threeparttable}
\usepackage{slashbox}
\usepackage{cite}
\renewcommand{\baselinestretch}{1.8}

\newtheorem{theorem}{Theorem}
\newtheorem{remark}{Remark}
\ifCLASSINFOpdf
\else
\fi
\begin{document}
%
\title{A dynamic state transition algorithm with application to sensor network localization}
%
%
%

\author{Xiaojun Zhou$^{1}$, Peng Shi$^{2}$, Cheng-Chew Lim$^{2}$, Chunhua Yang$^{1}$, Weihua Gui$^{1}$
\thanks{$^{1}$X.Zhou, C. Yang and W. Gui are with the School of Information Science and Engineering,
        Central South University, Changsha 410083, China.
       $^{2}$P. Shi and C.C. Lim are with the School of Electrical and Electronic Engineering, The University of Adelaide, Adelaide, SA 5005, Australia}
}

\maketitle

\begin{abstract}
The sensor network localization (SNL) problem is to reconstruct the positions of all the sensors in a network with
the given distance between pairs of sensors and within the radio range between them. It is proved that
the computational complexity of the SNL problem is NP-hard, and semi-definite programming or second-order cone programming relaxation methods are only able to solve some special problems of this kind.
In this study, a stochastic global optimization method called the state transition algorithm is introduced to solve the SNL problem without additional assumptions and conditions of the problem structure.
To transcend local optimality, a novel dynamic adjustment strategy called ``risk and restoration in probability"
is incorporated into the state transition algorithm. An empirical study is investigated to appropriately choose the ``risk probability" and ``restoration probability", yielding the dynamic state transition algorithm, which  is further improved by gradient-based refinement. The dynamic state transition algorithm  with refinement is applied to the SNL problem, and satisfactory experimental results have testified the effectiveness of the proposed approach.
\end{abstract}


%
\IEEEpeerreviewmaketitle

\section{Introduction}
%
%
%
%
The state transition algorithm has emerged in recent years as a novel stochastic method in global optimization, in which,
a solution to an optimization problem is regarded as a state,
and the update of current solution is regarded as a state transition \cite{zhou2012state,zhou2013comparative,zhou2014nonlinear}.
Using the state space representation, the state transition algorithm can describe solutions updating
in a unified framework, and the execution operators to update solutions are expressed as state transition matrices, which make it easy to understand and flexible to implement.
In the continuous state transition algorithm, there exist four state transformation operators, \textit{viz},
rotation, translation, expansion and axesion. These operators have special characteristics, covering the
local and global search capability. For example, the rotation operator has the ability to search in a hypersphere within a given radius, belonging to a local search; while the expansion operator has the probability to search in the whole space, belonging to a global search. The job specialization makes it convenient for users to manipulate
the state transition algorithm according to their demands. The strong global search ability and adaptability
of state transition algorithm have been demonstrated by comparison with other global optimization algorithms and several real-world applications \cite{zhou2012state,zhou2013comparative,zhou2014nonlinear}.

In recent decades, ad hoc wireless sensor networks have received considerable attention due to easy installation and simple operation \cite{la2013distributed,mao2007wireless,wang2009distributed,yoon2013efficient}. A typical sensor network consists of a large number of sensors, deployed in a geographical area. Sensor nodes
collect the local information (temperature, humidity, and vibration motion)  and communicate
with other neighboring nodes (two nodes are neighbors if the distance between them is below
some threshold, called the radio range). The sensor data from these nodes are
relevant only if we know what location they refer to; therefore, knowledge of the node positions
becomes imperative. Although locating these positions can be achieved
through manual configuration or by using the Global Positioning System (GPS) techniques, neither
methodology works well and both have physical limitations.
As a result, techniques to estimate node positions have been shifted to
develop methods that rely only on the measurements of distances
between neighboring nodes. The distance information could
be based on criteria such as time of arrival, time-difference
of arrival and received signal strength. Furthermore, it
is assumed that we already know the positions of a few sensor nodes (called anchors). For given positions of anchor nodes and relative distance between neighboring nodes,
the problem of finding the positions of all the sensor nodes is called the sensor network localization (SNL) problem \cite{liang2004gradient}.

The difficulty of locating the unknown sensors accurately is three-fold: 1) the distance
measurements may contain some noise or uncertainty; 2) it is not easy to identify the sufficient conditions for the sensor network to be localizable; and 3) the sensor network localization problem is proved to be NP-hard \cite{aspnes2004computational}.
The NP-hardness has led to efforts being directed at solving this problem approximately or solving it completely under certain conditions.
Semi-definite programming (SDP) relaxation has been widely used for the SNL problem \cite{biswas2004semidefinite,liang2004gradient,so2007theory,wang2008further}.
However, the solutions obtained by SDP relaxation are not generally optimal or are even infeasible; therefore, a
rounding technique is necessary to round the SDP solution to a suboptimal and feasible one.
Since the distance measurements inevitably contain noise or uncertainties, such methods to rounding the SDP solution may become not so robust and reliable.
Second-order cone programming (SOCP) relaxation has also found applications for the SNL problem \cite{srirangarajan2008distributed,tseng2007second}.
It is shown that even if it is weaker than the SDP relaxation, the SOCP relaxation has simpler structure and nicer properties that can make it useful as a problem preprocessor due to its faster speed.
Other gradient-based methods for the SNL problem can also be found in \cite{carter2006spaseloc,nie2009sum,ruan2014global,wu2012canonical,zhou2014canonical} and the references therein.

To the best our knowledge, there exist very few stochastic methods for the problem.
Some particle swarm optimization algorithms were used to solve the SNL problem \cite{gopakumar2009performance,kulkarni2011particle,low2008particle}, but the size of the problem is less than 100.
In this study, the state transition algorithm (STA) is applied to solve the SNL problem with the size more than 100 and is also scalable, this is one of the motivations of this paper.
Since STA is a stochastic global optimization algorithm, in principle, it has the capability to find a global solution to
the SNL problem without additional assumptions and conditions.
Nevertheless, as is known to us, it is inevitable to get trapped into local minima for most stochastic optimization algorithms.
To transcend local optimality, a dynamic adjustment strategy called ``risk and restoration in probability" is proposed to improve the global search ability. The ``risk in probability" means that a relatively worse solution
is accepted as next iterative with a probability $p_1$, while ``restoration in probability" means that the historical best solution is restored with another probability $p_2$. The values of these two probability are
investigated by an empirical study in this study.
is a good choice. With the proposed dynamic adjustment strategy, a dynamic state transition algorithm with refinement
is presented for the SNL problem.
Several experimental results have demonstrated the effectiveness of the proposed approach.

The main contribution of this paper is four-fold: 1) A fast rotation transformation operator is designed to reduce the computational complexity.
2) A dynamic state transition algorithm with ``risk and restoration in probability" strategy is proposed to escape from local optimality.
3) A good combination of the risk probability and the restoration probability is obtained by an empirical study.
4) The proposed dynamic state transition algorithm with refinement is successfully applied to the sensor network localization problem.

This paper is organized as follows. In section II, we give a brief review of the basic state transition algorithm (STA) and add some
modifications. Then, in section III, the local convergence analysis of the basic STA is discussed and a dynamic adjustment strategy is proposed to
improve its global search ability.  The proposed dynamic STA with refinement is applied to the sensor network localization
problem in section IV and the concluding
remarks are given in Section V.

\section{State Transition Algorithm}
In this paper, we focus on the following unconstrained optimization problem
\begin{equation}
\min_{\bm x \in \mathbb{R}^n} f(\bm x),
\end{equation}
where $\bm x$ is the decision variable vector, and $f$ is a single objective function.

In an iterative method, we update a current solution $\bm x_k$ to the next one $\bm x_{k+1}$. Similarly, in a state transition way, a solution can be regarded as a state, and the updating of a solution can be considered as a state transition process. \\
\indent On the basis of state space representation, the form of state transition algorithm can be described as follows:
\begin{equation}
\left \{ \begin{array}{ll}
\bm x_{k+1}= A_{k} \bm x_{k} + B_{k} \bm u_{k}\\
y_{k+1}= f(\bm x_{k+1})
\end{array} \right.,
\end{equation}
where $\bm x_{k}$ stands for a state, corresponding to a solution of the optimization problem; $\bm u_{k}$ is a function of $\bm x_{k}$ and history states; $y_{k}$ is the fitness value at $\bm x_k$;
$A_{k}$ and
$B_{k}$ are state transition matrices, which are usually some transformation operators;
 $f$ is the objective function or evaluation function.
\subsection{Original state transition algorithm}
\begin{figure}[!htbp]
  \centering
  \includegraphics[width=9cm]{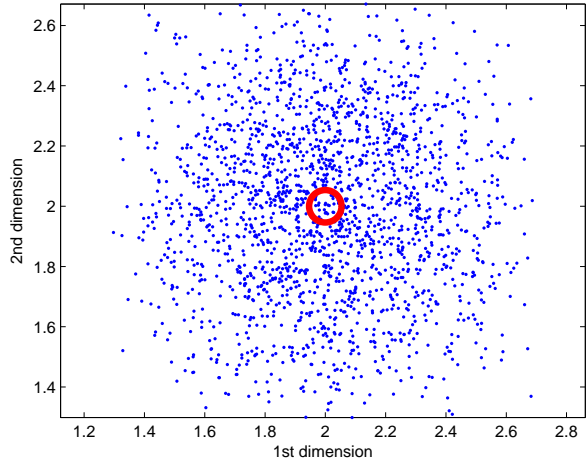}
  \caption{Illustration of the rotation transformation ($\alpha = 1$)}
  \label{rotation}
\end{figure}

\begin{figure}[!htbp]
  \centering
  \includegraphics[width=9cm]{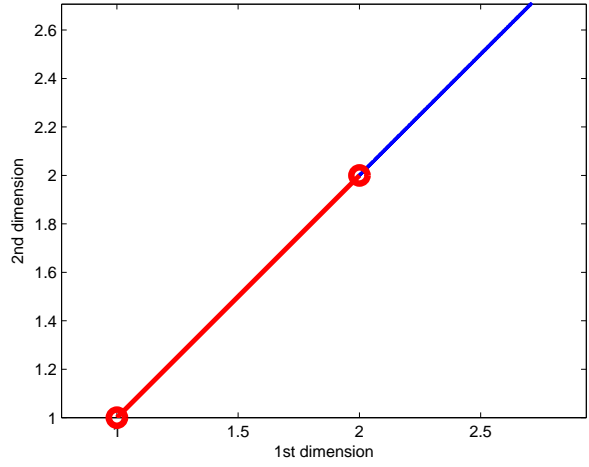}
  \caption{Illustration of the translation transformation ($\beta = 1$)}
  \label{translation}
\end{figure}

\begin{figure}[!htbp]
  \centering
  \includegraphics[width=9cm]{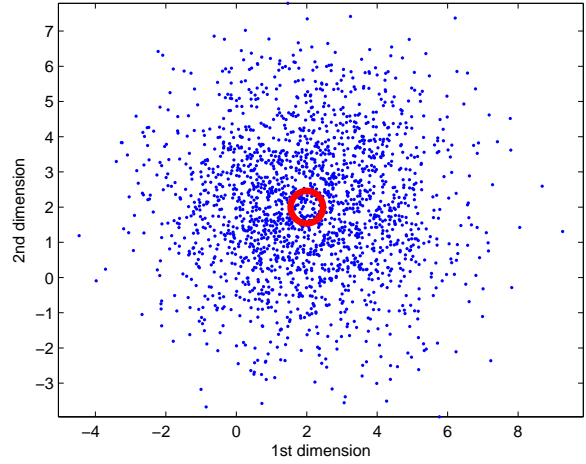}
  \caption{Illustration of the expansion transformation ($\gamma = 1$)}
  \label{expansion}
\end{figure}

\begin{figure}[!htbp]
  \centering
  \includegraphics[width=9cm]{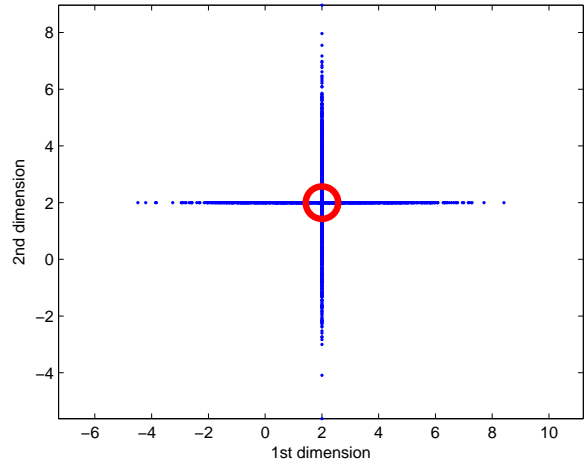}
  \caption{Illustration of the axesion transformation $(\delta = 1)$}
  \label{axesion}
\end{figure}

\indent Using space transformation, four special
state transformation operators are designed to solve the continuous function optimization problems.\\
(1) Rotation transformation
\begin{equation}
\bm x_{k+1}= \bm x_{k}+\alpha \frac{1}{n \|\bm x_{k}\|_{2}} R_{r} \bm x_{k},
\end{equation}
where $\alpha$ is a positive constant, called the rotation factor;
$R_{r}$ $\in$ $\mathbb{R}^{n\times n}$, is a random matrix with its entries being uniformly distributed random variables defined on the interval [-1, 1],
and $\|\cdot\|_{2}$ is the 2-norm of a vector. This rotation transformation
has the function of searching in a hypersphere with a maximal radius $\alpha$, as shown below
\begin{equation}
\begin{split}
\|\bm x_{k+1}-\bm x_{k}\|_{2} & = \|\alpha \frac{1}{n \|\bm x_{k}\|_{2}} R_{r} \bm x_{k} \|_{2}\\
                      & = \frac{\alpha}{n \|\bm x_{k}\|_{2}} \|R_{r} \bm x_{k}\|_{2}\\
                      & \leq \frac{\alpha}{n \|\bm x_{k}\|_{2}} \|R_{r}\|_{m_\infty}\|\bm x_{k}\|_{2} \leq \alpha.
\end{split}
\end{equation}
The illustration of
the rotation transformation in 2-D is given in Fig.\ref{rotation}, here, $\bm x_k = (2,2)$.
\\
(2) Translation transformation\\
\begin{equation}
\bm x_{k+1} = \bm x_{k}+  \beta  R_{t}  \frac{\bm x_{k}- \bm x_{k-1}}{\|\bm x_{k}- \bm x_{k-1}\|_{2}},
\end{equation}
where $\beta$ is a positive constant, called the translation factor; $R_{t}$ $\in \mathbb{R}$ is a uniformly distributed random variable defined on the interval [0,1]. Fig.\ref{translation} shows that
the translation transformation has the function of searching along a line from $\bm x_{k-1}(1,1)$ to $\bm x_{k}(2,2)$
 at the starting point $\bm x_{k}$ with a maximal length $\beta$.
\\
(3) Expansion transformation
\begin{equation}
\bm x_{k+1} = \bm x_{k}+  \gamma  R_{e} \bm x_{k},
\end{equation}
where $\gamma$ is a positive constant, called the expansion factor; $R_{e} \in \mathbb{R}^{n \times n}$ is a random diagonal
matrix with its entries obeying the Gaussian distribution. Fig.\ref{expansion} shows that the expansion transformation
has the function of expanding the entries in $\bm x_{k}$ to the range of [-$\infty$, +$\infty$], searching in the whole space.
\\
(4) Axesion transformation
\begin{equation}
\bm x_{k+1} = \bm x_{k}+  \delta  R_{a}  \bm x_{k}\\
\end{equation}
where $\delta$ is a positive constant, called the axesion factor; $R_{a}$ $\in \mathbb{R}^{n \times n}$ is a random diagonal matrix with its entries obeying the Gaussian distribution and only one random position having nonzero value. As illustrated in Fig.\ref{axesion}, the axesion transformation is aiming to search along the axes.

With those four state transformation operators, the original state transition algorithm can be described by the following pseudocode
\begin{algorithmic}[1]
\Repeat
    \If{$\alpha < \alpha_{\min}$}
    \State {$\alpha \gets \alpha_{\max}$}
    \EndIf
    \State {Best $\gets$ expansion(funfcn,Best,SE,$\beta$,$\gamma$)}
    \State {Best $\gets$ rotation(funfcn,Best,SE,$\alpha$,$\beta$)}
    \State {Best $\gets$ axesion(funfcn,Best,SE,$\beta$,$\delta$)}
    \State {$\alpha \gets \frac{\alpha}{\textit{fc}}$}
\Until{the specified termination criterion is met}
\end{algorithmic}
Note that $SE$ is called search enforcement, representing the times of transformation by a certain operator.

\begin{remark}
In the pseudocode of original state transition algorithm, the ``greedy criterion" is used to accept a new best solution from possible candidate states. The rotation factor $\alpha$ is decreasing from a maximum value to a minimum value in an
exponential way with the base $fc$, and other transformation factors are kept constant.
Due to its intrinsic properties, it can be seen that rotation is for exploitation (local search),
the expansion is for exploration (global search), the translation can be considered as line search which is applied only when a better solution is found, and the axesion is for strengthening the single dimensional search.
\end{remark}
\subsection{Fast rotation transformation}
\begin{figure}[!htbp]
  \centering
  \includegraphics[width=9cm]{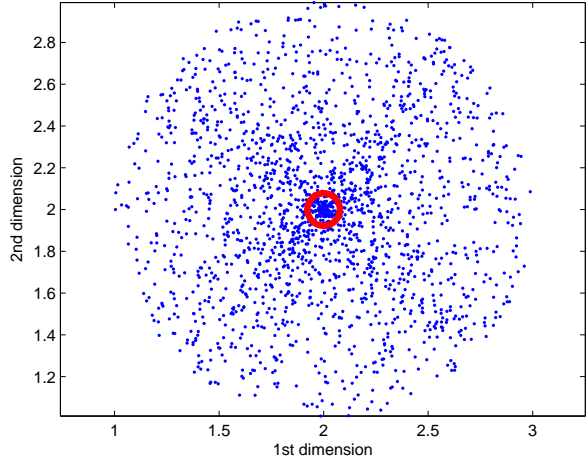}
  \caption{Illustration of the fast rotation transformation}
  \label{rotation_fast}
\end{figure}

\begin{figure}[!htbp]
  \centering
  \includegraphics[width=9cm]{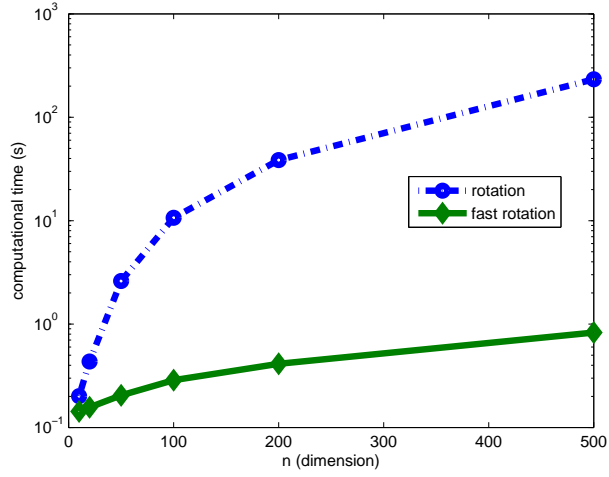}
  \caption{Comparison of the computational time for rotation operators}
  \label{comparison_rotation}
\end{figure}

To reduce the computational complexity, we propose a fast rotation transformation as
follows:

(5) Fast rotation transformation
\begin{equation}
\bm x_{k+1}= \bm x_{k} + \alpha \hat{R}_r \frac{\bm u}{\|\bm u\|_2},
\end{equation}
where $\hat{R}_r \in \mathbb{R}$ is a uniformly distributed random variable defined on the interval [-1,1],
and $\bm u \in \mathbb{R}^n$ is a vector with its entries being uniformly distributed random variables defined on the interval [-1,1]. It is easy to verify that
\begin{equation}
\begin{split}
\|\bm x_{k+1}-\bm x_{k}\|_{2} & = \Big\|\alpha \hat{R}_r \frac{\bm u}{\|\bm u\|_2} \Big\|_{2}\\
                      & = \alpha\|\hat{R}_r\|_{2}  \leq \alpha.
\end{split}
\end{equation}
The illustration of
the fast rotation transformation is given in Fig.\ref{rotation_fast}.
Compared with the original rotation operator,
the fast rotation transformation has low computational complexity since
the new random variable $\hat{R}_r$ is a scalar, not a matrix $R_r$ as in the previous rotation transformation.
Comparison of the computational time for two rotation operators can be found in
Fig.\ref{comparison_rotation}, where, at each fixed dimension, we run the two rotation operators for 10000 times on Intel(R) Core(TM) i3-2310M CPU @2.10GHz under Window 7 environment.
It is clear that the fast rotation operator can save significant amount of computing time.


The state transition algorithm (STA) with fast rotation transformation can be described by the following pseudocode
\begin{algorithmic}[1]
\Repeat
    \If{$\alpha(\beta,\gamma,\delta) < \alpha_{\min}(\beta_{\min},\gamma_{\min},\delta_{\min})$}
    \State {$\alpha(\beta,\gamma,\delta) \gets \alpha_{\max}(\beta_{\max},\gamma_{\max},\delta_{\max})$}
    \EndIf
    \State {Best $\gets$ expansion(funfcn,Best,SE,$\beta$,$\gamma$)}
    \State {Best $\gets$ rotation\_fast(funfcn,Best,SE,$\alpha$,$\beta$)}
    \State {Best $\gets$ axesion(funfcn,Best,SE,$\beta$,$\delta$)}
    \State {$\alpha(\beta,\gamma,\delta) \gets \frac{\alpha(\beta,\gamma,\delta)}{\textit{fc}}$}
\Until{the specified termination criterion is met}
\end{algorithmic}
\begin{remark}
We note that a fast rotation operator is introduced into the
new state transition transition. Furthermore, not only $\alpha$, but
also $\beta$, $\gamma$
and $\delta$ are
decreasing from a maximum value to a minimum value in an
exponential way with the base $fc$, which will be helpful for exploitation.
\end{remark}
\section{Dynamic state transition algorithm}
Although the new state transition algorithm
performs well for the majority of the benchmark functions, it behaves
weak for the Rosenbrock function, as shown in the experimental results of \cite{zhou2012state},
In fact, the Rosenbrock function is considered as a hard case for
most deterministic and stochastic optimization algorithms.
In this paper, we focus on how to jump out from local minima using a dynamic adjustment strategy, called ``risk and restoration in probability". The dynamic state transition algorithm risks accepting a relatively worse solution with a probability and restores the historical best solution with another probability.

\subsection{Convergence analysis of state transition algorithm}
Before presenting the dynamic state transition algorithm, we will
show that the previous proposed state transition algorithm
can at least converge to a local minimum.

From a control theory perspective, the evolution of the incumbent best state $\bm x_{k}^{*}$ in previous state transition algorithm using the ``greedy criterion" can be regarded as a discrete-time switched system
\begin{equation}
\bm x_{k+1}^{*} =
\left \{ \begin{array}{ll}
A_{k} \bm x_{k}^{*} + B_{k} \bm u_{k}, \mathrm{if}\; f(A_{k} \bm x_{k}^{*} + B_{k} \bm u_{k}) < f(\bm x_k^{*})\\
\bm x_{k}^{*}, \mathrm{otherwise}
\end{array} \right.
\end{equation}
A switched system contains some interesting stability phenomena. For instance, even when all the subsystems
are exponentially stable, a switched system may have divergent trajectories for certain switching signals.
That is to say, the stability of a switched system depends not only on the dynamics of each subsystem but also
on the properties of switching signals.
\begin{theorem}
The sequences $\{f(\bm x_k^{*})\}_{k=0}^{\infty}$ generated by the state transition algorithm can at least
converge to a local minimum, \textit{i.e.},
\begin{eqnarray}
\lim_{k \rightarrow \infty} f(\bm x_k^{*}) = f(\bar{\bm x}^{*})
\end{eqnarray}
where $\bar{\bm x}^{*}$ is a local minimum.
\end{theorem}
\proof
Since $f(\bm x_{k+1}^{*}) \leq f(\bm x_{k}^{*})$ and $f(\bm x_k^{*})$ is bounded below on $\mathbb{R}^n$, \textit{i.e.}, the global solution can be achieved, the sequence
$\{f(\bm x_k^{*})\}_{k=0}^{\infty}$ converges according to the \textit{monotone convergence theorem}.

Denote $\bar{\bm x}^{*}$ as the limiting point of the sequence $\{f(\bm x_k^{*})\}_{k=0}^{\infty}$.
Due to the rotation operator and accepting criteria used in the state transition algorithm, it is easy to find that
\begin{eqnarray}
f(\bar{\bm x}^{*}) \leq f(\bm x),\;\;\; \forall\; \|\bm x - \bar{\bm x}^{*}\| \leq \alpha_{\min}
\end{eqnarray}
that is to say, $\bar{\bm x}^{*}$ is a local minimum.
This completes the proof.
\subsection{Transcending local optimality}
Note that although the sequence
$\{f(\bm x_k^{*})\}_{k=0}^{\infty}$ converges when $k$ approaches infinity, it does not mean that the sequence converges to the global solution.
This phenomenon is called ``premature convergence" in the evolutionary computation community.
From a practical point of view,  ``premature convergence" is inevitable since there is no prior knowledge to judge whether a solution is indeed the global solution.
Thus, how to jump out of local minima becomes a significant issue.
``Risking a relatively bad solution in probability" is an effective
strategy to escape from local minima, as indicated in SA (simulated annealing) \cite{kirkpatrick1983optimization};
nevertheless, it is computationally expensive to achieve convergence \cite{mitra1985convergence}.

In this study, a new strategy called ``risk and restoration in probability" is proposed.
For each state transformation in the inner loop, a relatively worse solution is accepted based on
the ``risk probability". While in the outer loop, the best solution in external archive is restored to
update the incumbent current solution. The pseudocode for the dynamic state transition algorithm for
the outer loop is given as follows
\begin{algorithmic}[1]
\Repeat
    \State {[Best,fBest] $\gets$ expansion(funfcn,Best,SE,$\beta$,$\gamma$)}
    \State {[Best,fBest] $\gets$ rotation\_fast(funfcn,Best,SE,$\alpha$,$\beta$)}
    \State {[Best,fBest] $\gets$ axesion(funfcn,Best,SE,$\beta$,$\delta$)}
    \If{ fBest $<$ fBest$^{*}$}
    \State {Best$^{*}$ $\gets$ Best}
    \State {fBest$^{*}$ $\gets$ fBest}
    \EndIf
    \If{$rand < p_1$} \Comment{restoration in probability}
    \State {Best $\gets$ Best$^{*}$}
    \State {fBest $\gets$ fBest$^{*}$}
    \EndIf
\Until{the maximum number of iterations is met}
\end{algorithmic}
In this algorithm, Best and fBest are considered as the incumbent current solution and its function value,
Best$^{*}$ and fBest$^{*}$ are the best solution and its function value in history, and they are
kept in an external archive.

In the inner loop, taking expansion function for example, the pseudocode is given in the following
\begin{algorithmic}[1]
\State{State $\gets$ op\_expand(Best,SE,$\gamma$)}
\State{[newBest,fnewBest] $\gets$ fitness(funfcn,State)}
\If{fnewBest $<$ fBest}
    \State{Best $\gets$ newBest}
    \State{fBest $\gets$ fnewBest}
\Else
    \If{$rand < p_2$}\Comment{risk in probability}
    \State{Best $\gets$ newBest}
    \State{fBest $\gets$ fnewBest}
    \EndIf
\EndIf
\end{algorithmic}
\begin{remark}
In the above pseudocodes, $p_1$ and $p_2$ are the restoration probability and the risk probability, respectively.
Compared with the accepting criteria in simulating annealing, the
``risk and restoration in probability" strategy is easier to implement and has better convergent performance
since the best solution in history is always kept and restored frequently.
\end{remark}

\subsection{Benchmark functions test}
When evaluating the performance of the dynamic state transition algorithm,
some well-known benchmark functions are often used.
Listed below are the functions used in this work.

Spherical function
\begin{equation*}
f_1= \sum_{i=1}^n x_{i}^{2}, \; \; x_{i} \in [0,100],
\end{equation*}

Rastrigin function
\begin{equation*}
f_2=  \sum_{i=1}^n(x_{i}^{2}-10\cos(2 \pi x_{i})+10),\;\; x_{i} \in [0,5.12],
\end{equation*}

Griewank function
\begin{equation*}
f_3=\frac{1}{4000} \sum_{i=1}^n x_{i}^{2}-  \prod_{i}^{n} \cos|\frac{x_{i}}{\sqrt{i}}| + 1,\;\;  x_{i} \in [0,600],
\end{equation*}

Rosenbrock function
\begin{equation*}
f_4= \sum_{i=1}^n (100(x_{i+1}-x_{i}^2)^2 + (x_{i}-1)^2), \;\;  x_{i} \in [0,30],
\end{equation*}

Ackley function
\[
\begin{array}{l}
f_{5}=-20\exp(-0.2\sqrt{\frac{1}{n} \sum_{i=1}^n x_{i}^2})
\!-\!\exp(\frac{1}{n} \sum_{i=1}^n \cos(2\pi x_{i}))
+20+e, \;\; x_{i} \in [0,32].
\end{array}
\]

To investigate the effect of parameters in dynamical adjustment strategy of state transition algorithm,
some empirical study is arranged in order to find a satisfactory combination of the ``restoration probability" $p_1$ and the ``risk probability" $p_2$.

\begin{figure}[!htbp]
  \centering
  \includegraphics[width=8cm]{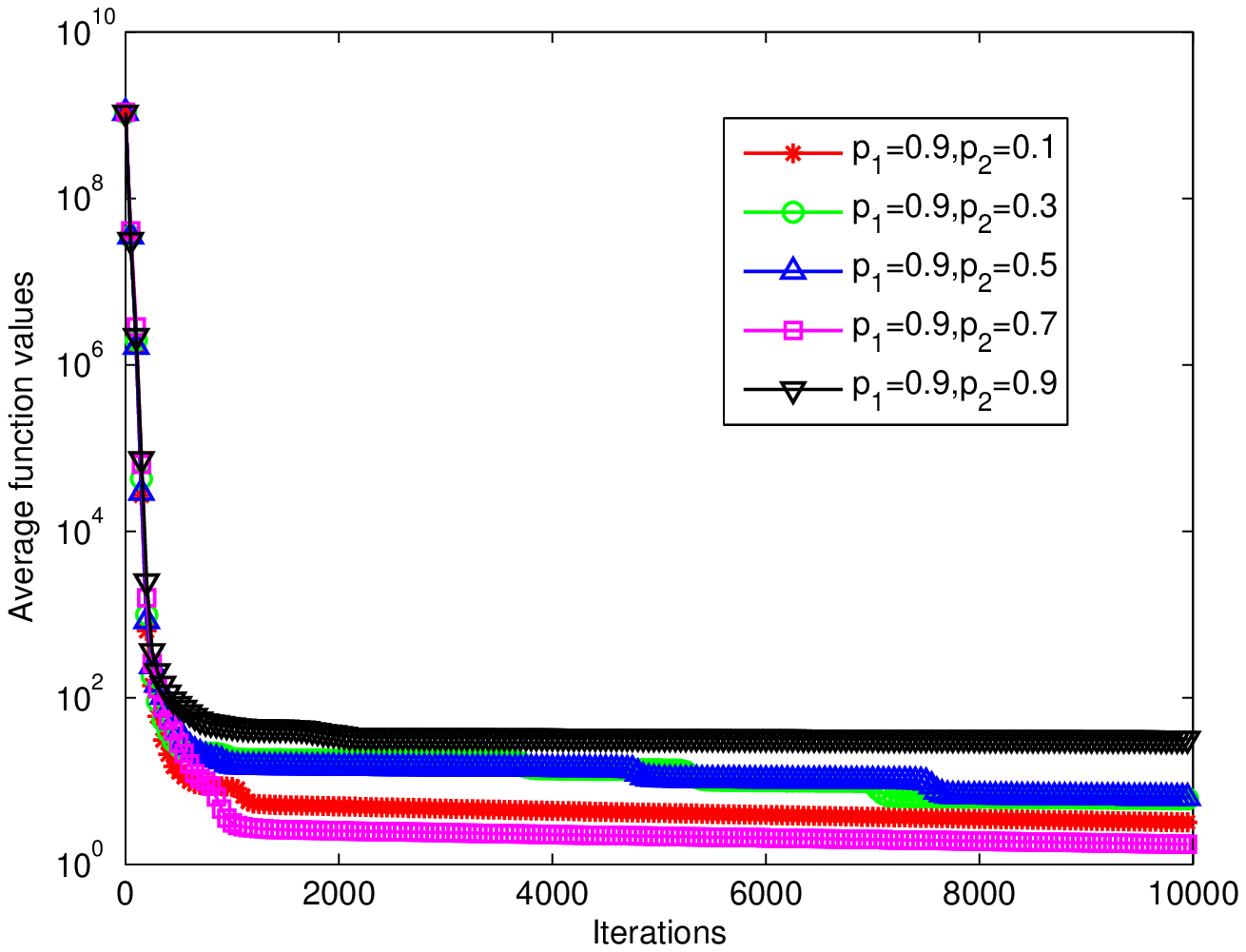}
  \includegraphics[width=8cm]{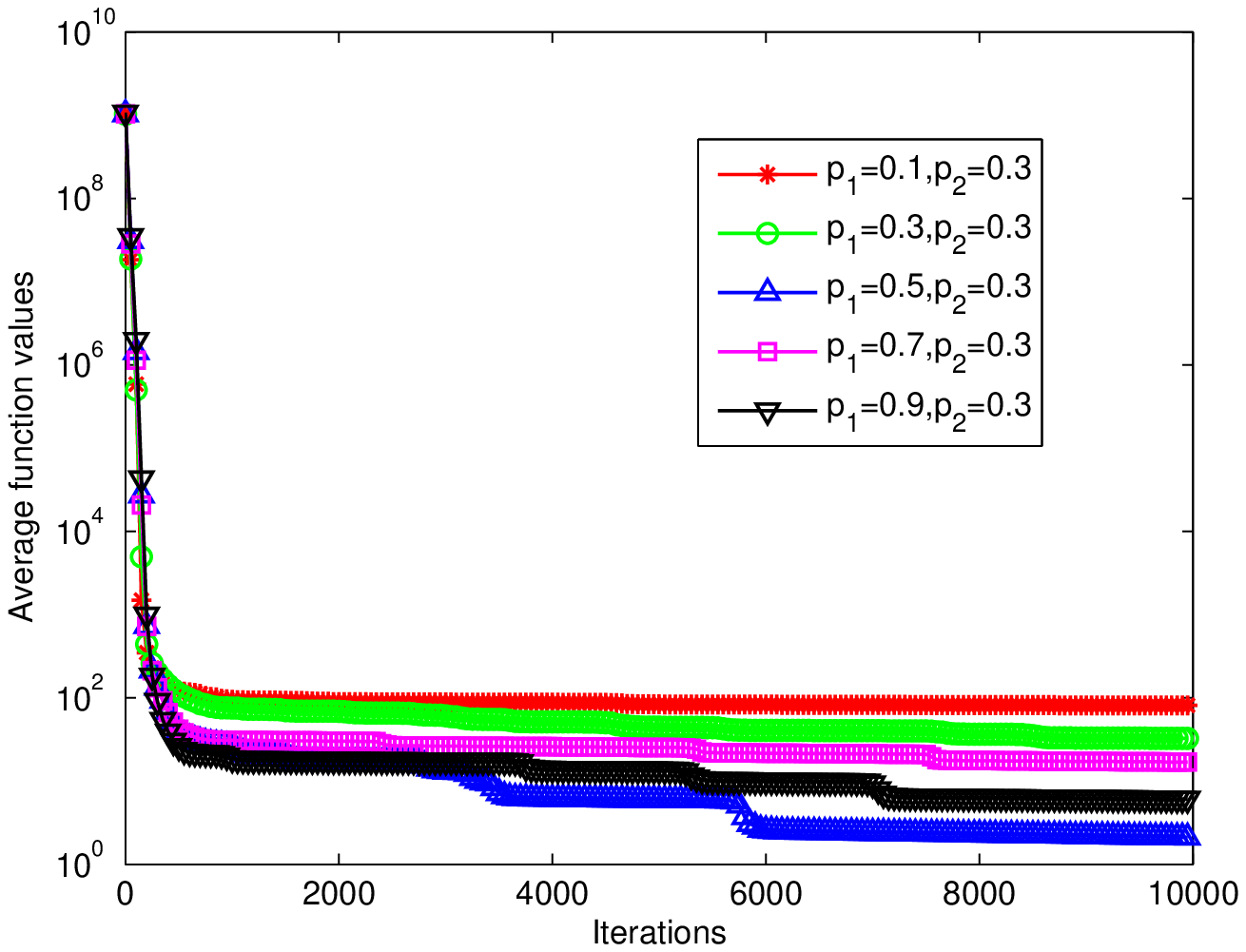}
  \caption{Iterative curves of the dynamic state transition algorithm with fixed ``risk probability" $p_1 = 0.9$ or ``restoration probability" $p_2 = 0.3$}
  \label{fig_dynamicparameters}
\end{figure}

First, with fixed dimension (problem size) and maximum number of iterations, some experimental tests
are performed on the benchmark functions with different groups of $(p_1,p_2)$.
From the experimental results, compared with the original STA (listed below in Table \ref{DSTA_benchmarks}),
we observe that for most groups of $(p_1,p_2)$, excluding $(p_1,p_2) = (0.1,0.9)$ and $(p_1,p_2) = (0.3,0.9)$, for all the benchmark functions except the Rosenbrock function,
the statistical performance of the dynamic STA and the original STA is almost the same when using the Wilcoxon rank sum test. However, for the Rosenbrock problem, the results of the dynamic STA with parameter groups $(p_1,p_2) = (0.3,0.1)$,  $(p_1,p_2) = (0.5,0.1)$,
 $(p_1,p_2) = (0.7,0.1)$,  $(p_1,p_2) = (0.9,0.1)$, $(p_1,p_2) = (0.5,0.3)$, $(p_1,p_2) = (0.9,0.3)$,
$(p_1,p_2) = (0.9,0.5)$ and $(p_1,p_2) = (0.9,0.7)$ are better or similar to those obtained by the original STA.
For a fixed ``restoration probability" $p_1 = 0.9$ or a fixed ``risk probability" $p_2 = 0.3$, the iterative curves of the dynamic state algorithm with the setting of other parameters are given in Fig.\ref{fig_dynamicparameters}. It can be seen that, either ``risk probability" or
``restoration probability" can influence the performance of the dynamic state transition algorithm.

\begin{figure}[!htbp]
  \centering
  \includegraphics[width=10cm]{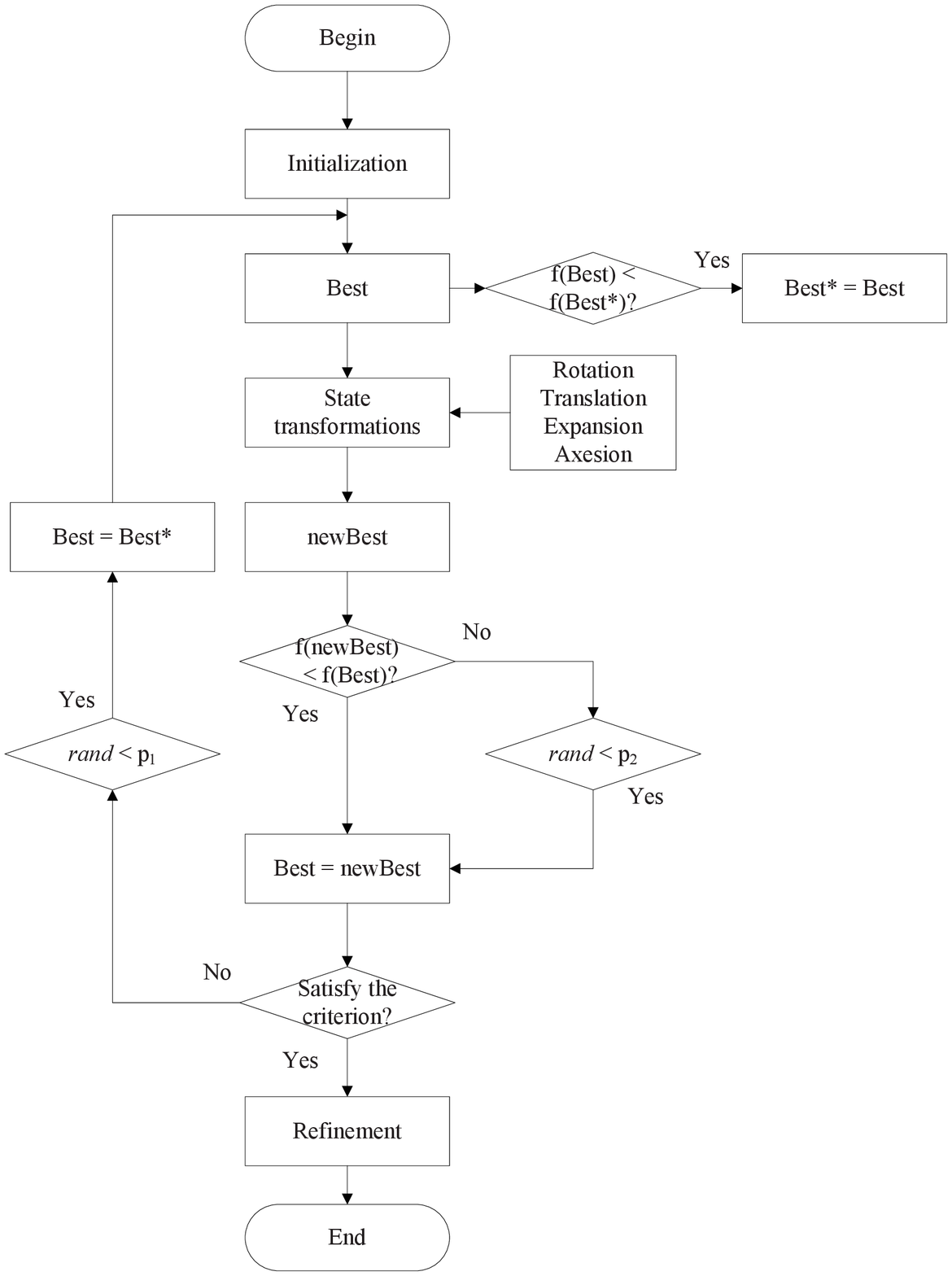}
  \caption{The flowchart of the dynamic state transition algorithm with refinement}
  \label{fig_flowchart}
\end{figure}

To further investigate the effect of $(p_1,p_2)$ on the Rosenbrock function, which belongs to the sum-of-squares problem (the sensor network localization problem studied latter also belongs to this kind), we increase the maximum number of iterations (from 100 to 500 times the dimension) and refine the solutions obtained by incorporating a gradient-based method.  As shown in Table \ref{DSTA_Rosenbrock_before} and Table \ref{DSTA_Rosenbrock_after}, the majority of the above mentioned parameter groups still perform well for the dynamic STA. The results after refinement are much closer to the true global solutions of the Rosenbrock problem, and the larger the maximum number of iterations, the closer to the true global solutions.
After taking stability and reliability into consideration, the parameter group $(p_1,p_2) = (0.9,0.3)$ is selected as an empirical best choice.
\begin{remark}
Through the benchmark function tests, the ``restoration probability" $p_1$ and the ``risk probability" $p_2$ in the dynamic state transition algorithm are determined empirically. The parameter group $(p_1,p_2) = (0.9,0.3)$ is also considered reasonable -- it tells us that a risk should be taken at a low probability, and the historical best should be restored at a high probability.

With a gradient-based technique incorporated for refinement, the flowchart of the dynamic state transition algorithm with refinement is illustrated in Fig.\ref{fig_flowchart}.
\end{remark}

\begingroup
\renewcommand{\baselinestretch}{1}%
\begin{table*}[!htbp]
\centering
\begin{threeparttable}[b]
\renewcommand{\arraystretch}{1.3}
\caption{Parametric study of the dynamic STA (dimension = 100, iterations = 10000)}
\scriptsize
\label{DSTA_benchmarks}
\begin{tabular}{{p{1.2cm}|p{0.8cm}|p{2.5cm}|p{2.5cm}|p{2.5cm}|p{2.7cm}|p{3.3cm}}}
\hline
Functions & \backslashbox[10mm]{$p_1$\kern-1em}{\kern-1em$p_2$}  & 0.1 & 0.3 & 0.5 & 0.7 & 0.9 \\
\hline
Spherical & 0.1 &  0 $\pm$ 0 $\approx$\tnote{1}  & 0 $\pm$ 0  $\approx$ & 0 $\pm$ 0 $\approx$ & 0 $\pm$ 0 $\approx$ & 1.6231e-15 $\pm$ 5.1263e-15 $-$  \\
          & 0.3 &  0 $\pm$ 0 $\approx$  & 0 $\pm$ 0 $\approx$ & 0 $\pm$ 0 $\approx$ & 0 $\pm$ 0 $\approx$ & 0 $\pm$ 0 $\approx$ \\
          & 0.5 &  0 $\pm$ 0 $\approx$  & 0 $\pm$ 0 $\approx$ & 0 $\pm$ 0 $\approx$ & 0 $\pm$ 0 $\approx$ & 0 $\pm$ 0 $\approx$ \\
          & 0.7 &  0 $\pm$ 0 $\approx$  & 0 $\pm$ 0 $\approx$ & 0 $\pm$ 0 $\approx$ & 0 $\pm$ 0 $\approx$ & 0 $\pm$ 0 $\approx$ \\
          & 0.9 &  0 $\pm$ 0 $\approx$  & 0 $\pm$ 0 $\approx$ & 0 $\pm$ 0 $\approx$ & 0 $\pm$ 0 $\approx$ & 0 $\pm$ 0 $\approx$ \\
\hline
Rastrigin & 0.1 &  0 $\pm$ 0 $\approx$  & 0 $\pm$ 0 $\approx$ & 0 $\pm$ 0 $\approx$ & 0 $\pm$ 0 $\approx$ & 3.3651e-12 $\pm$ 6.5316e-12 $-$ \\
          & 0.3 &  0 $\pm$ 0 $\approx$  & 0 $\pm$ 0 $\approx$ & 0 $\pm$ 0 $\approx$ & 0 $\pm$ 0 $\approx$ & 5.6843e-15 $\pm$ 2.5421e-14 $-$ \\
          & 0.5 &  0 $\pm$ 0 $\approx$  & 0 $\pm$ 0 $\approx$ & 0 $\pm$ 0 $\approx$ & 0 $\pm$ 0 $\approx$ & 0 $\pm$ 0 $\approx$ \\
          & 0.7 &  0 $\pm$ 0 $\approx$  & 0 $\pm$ 0 $\approx$ & 0 $\pm$ 0 $\approx$ & 0 $\pm$ 0 $\approx$ & 0 $\pm$ 0 $\approx$ \\
          & 0.9 &  0 $\pm$ 0 $\approx$  & 0 $\pm$ 0 $\approx$ & 0 $\pm$ 0 $\approx$ & 0 $\pm$ 0 $\approx$ & 1.1369e-14 $\pm$ 5.0842e-14 $-$\\
\hline
Griewank  & 0.1 & 0 $\pm$ 0  $\approx$  & 0 $\pm$ 0 $\approx$& 0 $\pm$ 0 $\approx$& 0 $\pm$ 0 $\approx$&  5.4956e-16 $\pm$ 1.3119e-15 $-$ \\
          & 0.3 & 0 $\pm$ 0 $\approx$   & 0 $\pm$ 0 $\approx$& 0 $\pm$ 0 $\approx$& 0 $\pm$ 0 $\approx$& 0 $\pm$ 0 $\approx$ \\
          & 0.5 & 0 $\pm$ 0 $\approx$   & 0 $\pm$ 0 $\approx$& 0 $\pm$ 0 $\approx$& 0 $\pm$ 0 $\approx$& 0 $\pm$ 0 $\approx$ \\
          & 0.7 & 0 $\pm$ 0 $\approx$   & 0 $\pm$ 0 $\approx$& 0 $\pm$ 0 $\approx$& 0 $\pm$ 0 $\approx$& 0 $\pm$ 0 $\approx$  \\
          & 0.9 & 0 $\pm$ 0 $\approx$   & 0 $\pm$ 0 $\approx$& 0 $\pm$ 0 $\approx$& 0 $\pm$ 0 $\approx$& 0 $\pm$ 0 $\approx$ \\
\hline
Rosenbrock& 0.1 & 34.9246 $\pm$ 48.8835 $-$  & 80.9344 $\pm$ 29.7257 $-$ &  83.5215 $\pm$ 19.7759 $-$ & 87.8305 $\pm$ 0.3616 $-$ & 88.0899 $\pm$ 0.2990 $-$  \\
          & 0.3 & 6.7223 $\pm$ 16.3449 $\approx$   & 32.5569 $\pm$ 45.5968 $-$ & 94.9096 $\pm$ 43.6650 $-$ & 111.6440 $\pm$ 56.4694 $-$ & 106.1008 $\pm$ 49.9611 $-$  \\
          & 0.5 & 6.2421 $\pm$ 17.6002 $\approx$   & 2.1033 $\pm$ 2.6610 $+$ & 28.3665 $\pm$ 54.6751 $-$ & 95.1745 $\pm$ 73.8811 $-$ & 91.6837 $\pm$ 77.8627 $-$ \\
          & 0.7 & 6.7115 $\pm$ 15.0144 $\approx$  & 16.8595 $\pm$ 34.9808 $-$ & 16.2994 $\pm$ 36.4491 $-$ & 24.7477 $\pm$ 35.8102 $-$ & 20.8548 $\pm$ 51.9346 $-$ \\
          & 0.9 &  3.1765 $\pm$ 3.0737 $+$    & 6.0835 $\pm$ 18.3434 $\approx$ & 6.2828 $\pm$ 18.1138 $\approx$ & 1.7381 $\pm$ 2.1587 $+$ & 31.9860 $\pm$ 61.3152 $-$ \\
\hline
Ackley    & 0.1 &  -8.8818e-16 $\pm$ 0  $\approx$ &  -8.8818e-16 $\pm$ 0 $\approx$ &  -8.8818e-16 $\pm$ 0 $\approx$ & -8.8818e-16 $\pm$ 0 $\approx$ & 2.8025e-9 $\pm$ 4.4539e-9 $-$ \\
          & 0.3 & -8.8818e-16 $\pm$ 0    $\approx$ & -8.8818e-16 $\pm$ 0 $\approx$ &  -8.8818e-16 $\pm$ 0 $\approx$ & -8.8818e-16 $\pm$ 0 $\approx$ & -8.8818e-16 $\pm$ 0 $\approx$ \\
          & 0.5 & -8.8818e-16 $\pm$ 0  $\approx$  & -8.8818e-16 $\pm$ 0 $\approx$ &  -8.8818e-16 $\pm$ 0 $\approx$ & -8.8818e-16 $\pm$ 0 $\approx$ & -8.8818e-16 $\pm$ 0 $\approx$ \\
          & 0.7 & -8.8818e-16 $\pm$ 0    $\approx$ & -8.8818e-16 $\pm$ 0 $\approx$ &  -8.8818e-16 $\pm$ 0 $\approx$ & -8.8818e-16 $\pm$ 0 $\approx$ &  -8.8818e-16 $\pm$ 0 $\approx$ \\
          & 0.9 & -8.8818e-16 $\pm$ 0 $\approx$   & -8.8818e-16 $\pm$ 0 $\approx$ &  -8.8818e-16 $\pm$ 0 $\approx$ & -8.8818e-16 $\pm$ 0 $\approx$&  -8.8818e-16 $\pm$ 0 $\approx$ \\
\hline
\hline
Methods&  & Spherical & Rastrigin & Griewank & Rosenbrock & Ackley\\
\hline
STA  & & 0 $\pm$ 0 \tnote{2} & 0 $\pm$ 0 & 0 $\pm$ 0 & 6.7952 $\pm$ 7.2544 & -8.8818e-16 $\pm$ 0\\
\hline
\end{tabular}
\begin{tablenotes}
\item [1] $+$, $-$ and $\approx$ denote that the performance of corresponding algorithm is better than, worse than, and similar to that of the STA by the Wilcoxon rank sum test
\item [2] $\pm$ denotes ``mean $\pm$ standard deviation"
\end{tablenotes}
\end{threeparttable}
\end{table*}
\endgroup

\begingroup
\renewcommand{\baselinestretch}{1}%
\begin{table*}[!htbp]
\centering
\renewcommand{\arraystretch}{1.1}
\caption{Performance of the dynamic STA for Rosenbrock function before refinement (dimension = 100)}
\scriptsize
\label{DSTA_Rosenbrock_before}
\newsavebox{\tablebox}
\begin{lrbox}{\tablebox}
\begin{tabular}{{p{1.25cm}|p{0.85cm}|p{3.0cm}|p{3.0cm}|p{3.0cm}|p{3.0cm}|p{3.0cm}}}
\hline
Iterations & \backslashbox[10mm]{$p_1$\kern-1em}{\kern-1em$p_2$}  & 0.1 & 0.3 & 0.5 & 0.7 & 0.9 \\
\hline
10000& 0.1 & 34.9246 $\pm$ 48.8835 $-$  & 80.9344 $\pm$ 29.7257 $-$ &  83.5215 $\pm$ 19.7759 $-$ & 87.8305 $\pm$ 0.3616 $-$ & 88.0899 $\pm$ 0.2990 $-$  \\
          & 0.3 & 6.7223 $\pm$ 16.3449 $\approx$   & 32.5569 $\pm$ 45.5968 $-$ & 94.9096 $\pm$ 43.6650 $-$ & 111.6440 $\pm$ 56.4694 $-$ & 106.1008 $\pm$ 49.9611 $-$  \\
          & 0.5 & 6.2421 $\pm$ 17.6002 $\approx$   & 2.1033 $\pm$ 2.6610 $+$ & 28.3665 $\pm$ 54.6751 $-$ & 95.1745 $\pm$ 73.8811 $-$ & 91.6837 $\pm$ 77.8627 $-$ \\
          & 0.7 & 6.7115 $\pm$ 15.0144 $\approx$  & 16.8595 $\pm$ 34.9808 $-$ & 16.2994 $\pm$ 36.4491 $-$ & 24.7477 $\pm$ 35.8102 $-$ & 20.8548 $\pm$ 51.9346 $-$ \\
          & 0.9 &  3.1765 $\pm$ 3.0737 $+$    & 6.0835 $\pm$ 18.3434 $\approx$ & 6.2828 $\pm$ 18.1138 $\approx$ & 1.7381 $\pm$ 2.1587 $+$ & 31.9860 $\pm$ 61.3152 $-$ \\
\hline
 20000& 0.1 &  36.6229 $\pm$ 41.1742 $-$  & 70.7605 $\pm$ 30.1326 $-$ & 80.3873 $\pm$ 18.9506 $-$ & 84.8749 $\pm$ 0.2333 $-$ & 85.4611 $\pm$ 0.2593 $-$ \\
          & 0.3 &  1.3527$\pm$ 1.6541 $+$  & 22.8921 $\pm$ 41.2240 $-$ & 62.3492 $\pm$ 58.6085 $-$ & 84.2958 $\pm$ 25.5081 $-$ & 100.3479 $\pm$ 44.1913 $-$ \\
          & 0.5 &  6.7016 $\pm$ 18.1658 $-$  & 1.3524 $\pm$ 1.6784 $+$ & 21.9584 $\pm$ 47.6510 $-$ & 48.4355 $\pm$ 58.3758 $-$ & 69.9912 $\pm$ 63.0667 $-$ \\
          & 0.7 &  1.8801 $\pm$ 2.0266 $+$  & 6.4236 $\pm$ 19.8274 $-$ & 6.8750 $\pm$ 21.5073 $-$ & 2.2362 $\pm$ 1.8067 $+$ & 21.3934 $\pm$ 42.2174 $-$ \\
          & 0.9 &  2.1100 $\pm$ 1.5590 $+$  & 1.9741 $\pm$ 1.9684 $+$ & 2.2778 $\pm$ 1.8401 $+$ & 11.8520 $\pm$ 34.9652 $-$ & 18.9113 $\pm$ 38.2404 $-$ \\
\hline
30000  & 0.1 & 8.1871 $\pm$ 22.3389  $-$  & 59.4579 $\pm$ 34.4171 $-$ & 81.9436 $\pm$ 0.8132 $-$ & 82.8267 $\pm$ 0.2116 $-$ &  83.5295 $\pm$ 0.2513 $-$ \\
          & 0.3 & 2.4303 $\pm$ 3.5954 $\approx$   & 30.1052 $\pm$ 44.2706 $-$ & 59.8675 $\pm$ 39.3182 $-$ & 83.9215 $\pm$ 39.0975 $-$ & 84.6498 $\pm$ 14.2251 $-$ \\
          & 0.5 & 1.0879 $\pm$ 1.3156 $+$   & 8.7277 $\pm$ 20.7838 $-$ & 18.0923 $\pm$ 41.2037 $-$ & 61.8106 $\pm$ 70.4271 $-$ & 74.3550 $\pm$ 64.2713 $-$ \\
          & 0.7 & 0.7939 $\pm$ 0.8984 $+$   & 1.3687 $\pm$ 1.1653 $+$& 9.9566 $\pm$ 24.4415 $-$ & 10.3543 $\pm$ 27.5641 $-$ & 33.5537 $\pm$ 49.5459 $-$ \\
          & 0.9 & 1.2680 $\pm$ 1.1676 $+$   & 1.2800 $\pm$ 1.1730 $+$ & 1.0241 $\pm$ 1.1646 $+$& 0.7466 $\pm$ 0.9355 $+$& 1.3722 $\pm$ 1.0550 $+$ \\
\hline
40000& 0.1 & 15.5817 $\pm$ 27.3889 $-$  & 57.2248 $\pm$ 33.4997 $-$ &  75.8210 $\pm$ 17.6466 $-$ & 81.0726 $\pm$ 0.3210 $-$ & 81.9298 $\pm$ 0.2200 $-$  \\
          & 0.3 & 0.6875 $\pm$ 0.8192 $+$   & 9.3453 $\pm$ 23.8299 $-$ & 46.0035 $\pm$ 46.7455 $-$ & 69.6654 $\pm$ 39.9633 $-$ & 84.8464 $\pm$ 21.9926 $-$  \\
          & 0.5 & 1.0538 $\pm$ 1.1154 $+$   & 13.0401 $\pm$ 29.9575 $-$ & 11.1267 $\pm$ 27.4796 $-$ & 17.6139 $\pm$ 30.8468 $-$ & 46.1817 $\pm$ 47.4052 $-$ \\
          & 0.7 & 0.9147 $\pm$ 0.6730 $+$  & 1.3028 $\pm$ 0.9836 $\approx$ & 0.5366 $\pm$ 0.7615 $+$ & 0.7609 $\pm$ 0.8081 $+$ & 27.7787 $\pm$ 46.9281 $-$ \\
          & 0.9 &  1.0351 $\pm$ 0.8745 $+$    & 0.7861$\pm$ 0.7066 $+$ & 0.5277 $\pm$ 0.8042 $+$ & 1.2816 $\pm$ 0.7297 $\approx$ & 10.2427 $\pm$ 29.4908 $-$ \\
\hline
50000     & 0.1 &  3.9619 $\pm$ 14.1378  $-$ &  64.2790 $\pm$ 27.5534 $-$ &  77.8868 $\pm$ 0.3830 $-$ & 79.5461 $\pm$ 0.2291 $-$ & 80.5820 $\pm$ 0.2589 $-$ \\
          & 0.3 & 0.7379 $\pm$ 0.7578    $+$ & 15.7587 $\pm$ 30.7724 $-$ &  26.7952 $\pm$ 36.4140 $-$ & 42.4749 $\pm$ 36.2685 $-$ & 79.6220 $\pm$ 37.3989 $-$ \\
          & 0.5 & 0.5307 $\pm$ 0.5997  $+$  & 0.5488 $\pm$ 0.7058 $+$ &  4.4704 $\pm$ 16.5813 $-$ & 26.8494 $\pm$ 42.6035 $-$ & 66.1161 $\pm$ 56.9610 $-$ \\
          & 0.7 & 0.6374 $\pm$ 0.5889    $+$ & 0.5150 $\pm$ 0.5831 $+$ &  5.5514 $\pm$ 14.3693 $-$ & 0.6901 $\pm$ 0.9297 $+$ &  10.1952 $\pm$ 25.8581 $-$ \\
          & 0.9 & 0.4749 $\pm$ 0.5146 $+$   & 0.6032 $\pm$ 0.6193 $+$ &  0.6189 $\pm$ 0.5609 $+$ & 1.3870 $\pm$ 3.8526 $-$&  1.2649 $\pm$ 2.1263 $-$ \\
\hline
\hline
Methods &   & 10000 & 20000 & 30000 & 40000 & 50000\\
\hline
STA  & & 6.7952 $\pm$ 7.2544 & 5.7962 $\pm$ 2.8367 & 2.8918 $\pm$ 1.9632 & 1.6408 $\pm$ 1.1822 & 1.0805 $\pm$ 0.7956\\
\hline
\end{tabular}
\end{lrbox}
\scalebox{0.87}{\usebox{\tablebox}}
\begin{tablenotes}
\item [] $+$, $-$ and $\approx$ denote that the performance of corresponding algorithm is better than, worse than, and similar to that of the STA by the Wilcoxon rank sum test
\item [] $\pm$ denotes ``mean $\pm$ standard deviation"
\end{tablenotes}
\end{table*}
\endgroup

\begingroup
\renewcommand{\baselinestretch}{1}%
\begin{table*}[!htbp]
\centering
\renewcommand{\arraystretch}{1.3}
\caption{Performance of the dynamic STA for Rosenbrock function after refinement (dimension = 100)}
\scriptsize
\label{DSTA_Rosenbrock_after}
\begin{lrbox}{\tablebox}
\begin{tabular}{{p{0.95cm}|p{0.85cm}|p{3.2045cm}|p{3.2045cm}|p{3.1cm}|p{3.08cm}|p{3.0cm}}}
\hline
Iterations & \backslashbox[8mm]{$p_1$\kern-1em}{\kern-1em$p_2$}  & 0.1 & 0.3 & 0.5 & 0.7 & 0.9 \\
\hline
10000& 0.1 & 21.1587 $\pm$ 36.1116 $-$  & 53.1184 $\pm$ 29.5350 $-$ &  63.3643 $\pm$ 20.1991 $-$ & 70.5142 $\pm$ 1.3252 $-$ & 70.2448 $\pm$ 1.0381 $-$  \\
          & 0.3 & 0.0204 $\pm$ 0.0912 $\approx$   & 7.5377 $\pm$ 22.0972 $-$ & 42.4304 $\pm$ 38.3322 $-$ & 50.5568 $\pm$ 42.2048 $-$ & 60.5727 $\pm$ 31.6022 $-$  \\
          & 0.5 & 0.2802 $\pm$ 1.2511 $-$   & 4.0550e-4 $\pm$ 0.0018 $+$ & 5.7047 $\pm$ 17.5013 $-$ & 51.1806 $\pm$ 54.7261 $-$ & 47.0276 $\pm$ 46.6980 $-$ \\
          & 0.7 & 5.7870e-4 $\pm$ 0.0026 $+$  & 2.6333e-6 $\pm$ 8.7040e-6 $+$ & 9.7799e-4 $\pm$ 0.0044 $+$ & 10.3570 $\pm$ 25.3780 $-$ & 10.3748 $\pm$ 25.3956 $-$ \\
          & 0.9 &  0.0013 $\pm$ 0.0057 $\approx$    & 5.7115e-5 $\pm$ 2.5542e-4 $+$ & 0.0920 $\pm$ 0.4115 $\approx$ & 9.3007e-10 $\pm$ 1.4912e-9 $+$ & 4.2359 $\pm$ 16.8367 $-$ \\
\hline
 20000& 0.1 &  17.6076 $\pm$ 31.7114 $-$  & 57.2231 $\pm$ 24.6991 $-$ & 64.7344 $\pm$ 15.3070 $-$ & 65.4911 $\pm$ 13.0773 $-$ & 68.9557 $\pm$ 1.0903 $-$ \\
          & 0.3 &  1.2964e-9 $\pm$ 1.7722e-9 $+$  & 4.5672e-9 $\pm$ 1.5513e-8 $-$ & 27.0481 $\pm$ 33.2740 $-$ & 53.1332 $\pm$ 26.1597 $-$ & 65.3754 $\pm$ 24.7380 $-$ \\
          & 0.5 &  3.3265 $\pm$ 14.8456 $-$  & 9.1372e-9 $\pm$ 3.6650e-8 $+$ & 5.4289 $\pm$ 18.2178 $-$ & 25.2178 $\pm$ 41.7371 $-$ & 28.3450 $\pm$ 36.0765 $-$ \\
          & 0.7 &  1.7181e-9 $\pm$ 2.9007e-9 $+$  & 0.2978 $\pm$ 1.3318 $-$ & 3.9112 $\pm$ 17.4915 $-$ & 2.3355e-9 $\pm$ 2.3355e-9 $+$ & 4.3373 $\pm$ 19.3952 $-$ \\
          & 0.9 &  2.5848e-9 $\pm$ 4.0031e-9 $+$  & 6.4586e-9 $\pm$ 2.4736e-8 $+$ & 9.7906e-10 $\pm$ 1.4118e-9 $+$ & 1.8054 $\pm$ 7.5612 $-$ & 7.2380e-8 $\pm$ 3.2018e-7 $+$ \\
\hline
30000  & 0.1 & 5.8828 $\pm$ 18.1071  $-$  & 45.2268 $\pm$ 30.3826 $-$ & 66.6607 $\pm$ 0.9523 $-$ & 66.9913 $\pm$ 0.9594 $-$ &  67.1216 $\pm$ 1.0454 $-$ \\
          & 0.3 & 0.1111 $\pm$ 0.4969 $-$   & 8.9327 $\pm$ 21.3730 $-$ & 28.6617 $\pm$ 36.1546 $-$ & 53.0224 $\pm$ 27.8523 $-$ & 61.5713 $\pm$ 18.5342 $-$ \\
          & 0.5 & 9.1503e-10 $\pm$ 1.4908e-9 $\approx$   & 2.8200 $\pm$ 12.6115 $-$ & 10.7349 $\pm$ 25.6296 $-$ & 21.6996 $\pm$ 32.4225 $-$ & 32.2821 $\pm$ 41.0405 $-$ \\
          & 0.7 & 4.9650e-10 $\pm$ 5.4828e-10 $\approx$   & 1.0867e-9 $\pm$ 1.4004e-9 $\approx$& 0.2918
 $\pm$ 1.3052 $-$ & 3.8141 $\pm$ 13.5439 $-$ & 7.9870 $\pm$ 22.5467 $-$ \\
          & 0.9 & 1.6998e-9 $\pm$ 2.7782e-9 $\approx$   & 1.4478e-9 $\pm$ 2.1626e-9 $\approx$ & 8.5037e-10 $\pm$ 1.5985e-9 $\approx$& 2.4864e-7 $\pm$ 1.1101e-6 $-$& 2.2585e-9
 $\pm$ 6.0799e-9 $\approx$ \\
\hline
40000& 0.1 & 10.7407 $\pm$ 21.8229 $-$  & 40.3449 $\pm$ 30.3853 $-$ &  61.1384 $\pm$ 14.4172 $-$ & 65.1246 $\pm$ 1.0444 $-$ & 65.3149 $\pm$ 0.9135 $-$  \\
          & 0.3 & 4.7625e-10 $\pm$ 5.8946e-10 $\approx$   & 6.6588 $\pm$ 18.9926 $-$ & 17.8645 $\pm$ 27.0329 $-$ & 40.9452 $\pm$ 29.2316 $-$ & 56.2615 $\pm$ 19.3012 $-$  \\
          & 0.5 & 1.2431e-9 $\pm$ 2.7342e-9 $\approx$   & 0.2716 $\pm$ 0.8660 $-$ & 4.8668 $\pm$ 19.8535 $-$ & 6.0495 $\pm$ 17.4021 $-$ & 27.8799 $\pm$ 33.3700 $-$ \\
          & 0.7 & 9.2246e-10 $\pm$ 1.4431e-9 $\approx$  & 1.4476e-9 $\pm$ 2.2476e-9 $\approx$ & 7.7317e-10 $\pm$ 1.6800e-9 $\approx$ & 1.0940e-9 $\pm$ 1.8583e-9 $\approx$ & 6.9124 $\pm$ 17.7629 $-$ \\
          & 0.9 & 7.0819e-10 $\pm$ 7.4544e-10 $\approx$    & 8.6857e-10 $\pm$ 1.2423e-9 $\approx$ & 8.6832e-10 $\pm$ 1.3687e-9 $\approx$ & 1.4373e-9 $\pm$ 1.8758e-9 $\approx$ & 1.3033 $\pm$ 5.8287 $-$ \\
\hline
50000     & 0.1 &  2.4878 $\pm$ 11.1259  $-$ &  42.2141 $\pm$ 28.4066 $-$ &  62.6582 $\pm$ 0.9551 $-$ & 64.0278 $\pm$ 0.7651 $-$ & 64.6633 $\pm$ 1.1308 $-$ \\
          & 0.3 & 7.3724e-10 $\pm$ 9.7763e-10    $\approx$ & 5.6969 $\pm$ 17.7139 $-$ &  14.1997 $\pm$ 24.6097 $-$ & 19.8874 $\pm$ 26.7939 $-$ & 54.9271 $\pm$ 20.6629 $-$ \\
          & 0.5 & 4.7795e-10 $\pm$ 5.8795e-10  $\approx$  & 1.7176e-9 $\pm$ 3.1498e-9 $\approx$ &  3.2398e-9 $\pm$ 1.2368e-8 $\approx$ & 8.5287 $\pm$ 20.8636 $-$ & 34.2486 $\pm$ 38.3122 $-$ \\
          & 0.7 & 2.0098e-9 $\pm$ 4.6715e-9    $\approx$ & 6.9608e-10 $\pm$ 6.2422e-10 $\approx$ &  3.0832 $\pm$ 10.7902 $-$ & 1.5723e-9 $\pm$ 3.4520e-9 $\approx$ &  4.5651 $\pm$ 18.0369 $-$ \\
          & 0.9 & 9.7722e-10 $\pm$ 1.5721e-9 $\approx$   & 8.5754e-10 $\pm$ 8.6167e-10 $\approx$ &  9.6368e-10 $\pm$ 1.4701e-9 $\approx$ & 0.1362 $\pm$ 0.6091 $-$&  1.5806e-9 $\pm$ 2.4191e-9 $\approx$ \\
\hline
\hline
Methods &   & 10000 & 20000 & 30000 & 40000 & 50000\\
\hline
STA  & & 0.0806 $\pm$ 0.2253 & 3.1866e-4 $\pm$ 0.0011 & 3.0420e-9 $\pm$ 4.8701e-9 & 2.7837e-9 $\pm$ 4.8575e-9 & 9.7776e-10 $\pm$ 1.2594e-9\\
\hline
\end{tabular}
\end{lrbox}
\scalebox{0.87}{\usebox{\tablebox}}
\begin{tablenotes}
\item [] $+$, $-$ and $\approx$ denote that the performance of corresponding algorithm is better than, worse than, and similar to that of the STA by the Wilcoxon rank sum test
\item [] $\pm$ denotes ``mean $\pm$ standard deviation"
\end{tablenotes}
\end{table*}
\endgroup

\section{Application for the sensor network localization}
In this section, the proposed dynamic state transition algorithm with refinement is applied to solve the sensor network localization problem.
The mathematical model of the sensor network localization problem can be described as follows.

Given $m$ anchor points $\bm a_1, \cdots, \bm a_m \in \mathbb{R}^d$ ($d$ is usually 2 or 3),
the distance $d_{ij}$ (Euclidean distance) between the $i$th and $j$th anchor points if $(i,j) \in N_x$, and the distance $e_{ik}$ between the $i$th sensor and $k$th anchor points
if $(i,k) \in N_a$, with $N_x = \{(i,j): \|\bm x_i - \bm x_j\| = d_{ij} \leq r_d\}$ and $N_a = \{(i,k): \|\bm x_i - \bm a_k\| = e_{ik} \leq r_d\}$, where, $r_d$ is the radio range.
The sensor network localization problem is to find $n$ distinct sensor points $\bm x_i, i = 1, \cdots, n$, such that
\begin{equation}
\label{eqwsn}
\|\bm x_i - \bm x_j\|^2 \!=\! d_{ij}^2, \forall (i,j) \in N_x, \nonumber \\
\|\bm x_i - \bm a_k\|^2 \!=\! e_{ik}^2, \forall (i,k) \in N_a.
\end{equation}

Usually, the distances $d_{ij}$ and $e_{ij}$ may contain noise, making equations (\ref{eqwsn}) infeasible; therefore, using the least squares method, we formulate the localization problem as
the following nonconvex optimization problem
\begin{equation}
\min \sum_{(i,j) \in N_x} (\|\bm x_i \!-\! \bm x_j\|^2 \!-\! d_{ij}^2)^2 \!+\! \sum_{(i,k) \in N_a} (\|\bm x_i \!-\! \bm a_k\|^2 \!-\! e_{ik}^2)^2.
\end{equation}

\begin{figure}[!htbp]
  \centering
  \includegraphics[width=8cm]{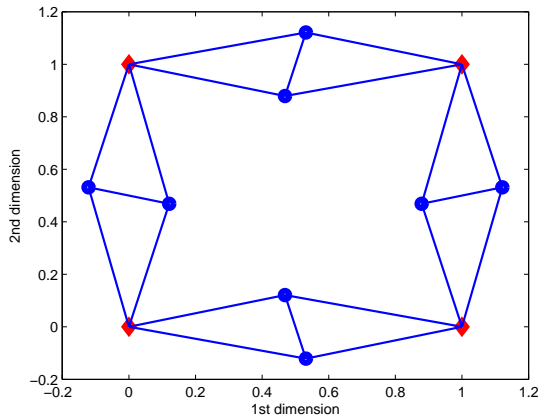}
  \caption{Network topology of the illustrative example}
  \label{fig_illustrative_example}
\end{figure}

We first consider the following illustrative example. The
network topology is given in Fig.\ref{fig_illustrative_example}, and
there are 4 anchors with known positions. We need to decide the positions of 8 sensors such that
the following distances between two sensors or between a sensor and an anchor are
satisfied:
\begin{eqnarray*}
&&\bm a_1 = (0,0), \bm a_2 = (0,1),  \bm a_3 = (1,0), \bm a_4 = (1,1), \\
&&d_{12} = 1/4,  d_{34} = 1/4,  d_{56} = 1/4,  d_{78} = 1/4, \\
&&e_{11} = \frac{\sqrt{15}}{8},  e_{13} = \frac{\sqrt{19}}{8},  e_{21} = \frac{\sqrt{19}}{8},  e_{23} = \frac{\sqrt{15}}{8}, \\
&&e_{33} = \frac{\sqrt{15}}{8},  e_{34} = \frac{\sqrt{19}}{8},  e_{43} = \frac{\sqrt{19}}{8},  e_{44} = \frac{\sqrt{15}}{8}, \\
&&e_{52} = \frac{\sqrt{19}}{8},  e_{54} = \frac{\sqrt{15}}{8},  e_{62} = \frac{\sqrt{15}}{8},  e_{64} = \frac{\sqrt{19}}{8}, \\
&&e_{71} = \frac{\sqrt{19}}{8},  e_{72} = \frac{\sqrt{15}}{8},  e_{81} = \frac{\sqrt{15}}{8},  e_{82} = \frac{\sqrt{19}}{8}. \\
\end{eqnarray*}

\begin{figure}[!htbp]
  \centering
  \includegraphics[width=8cm]{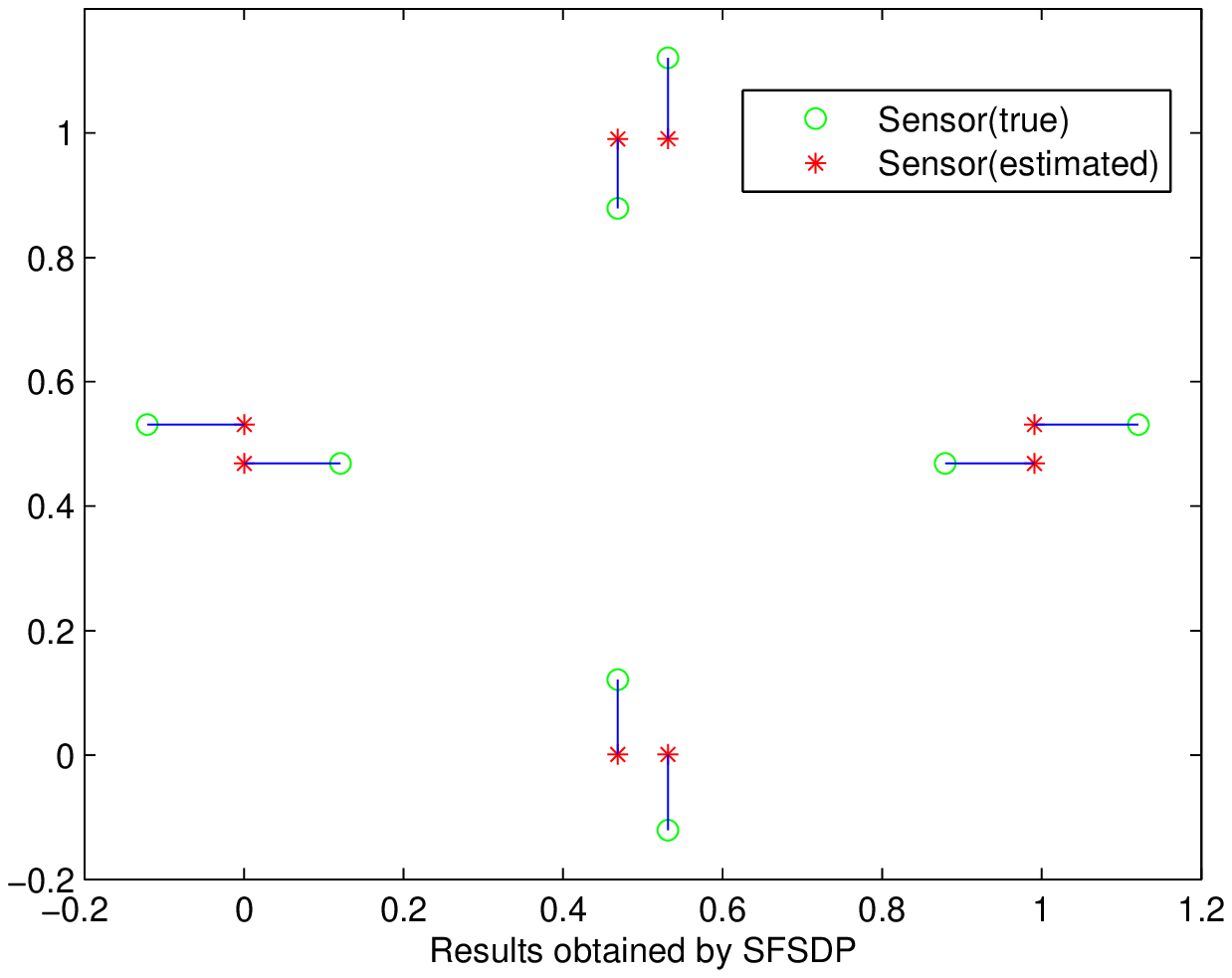}
  \includegraphics[width=8cm]{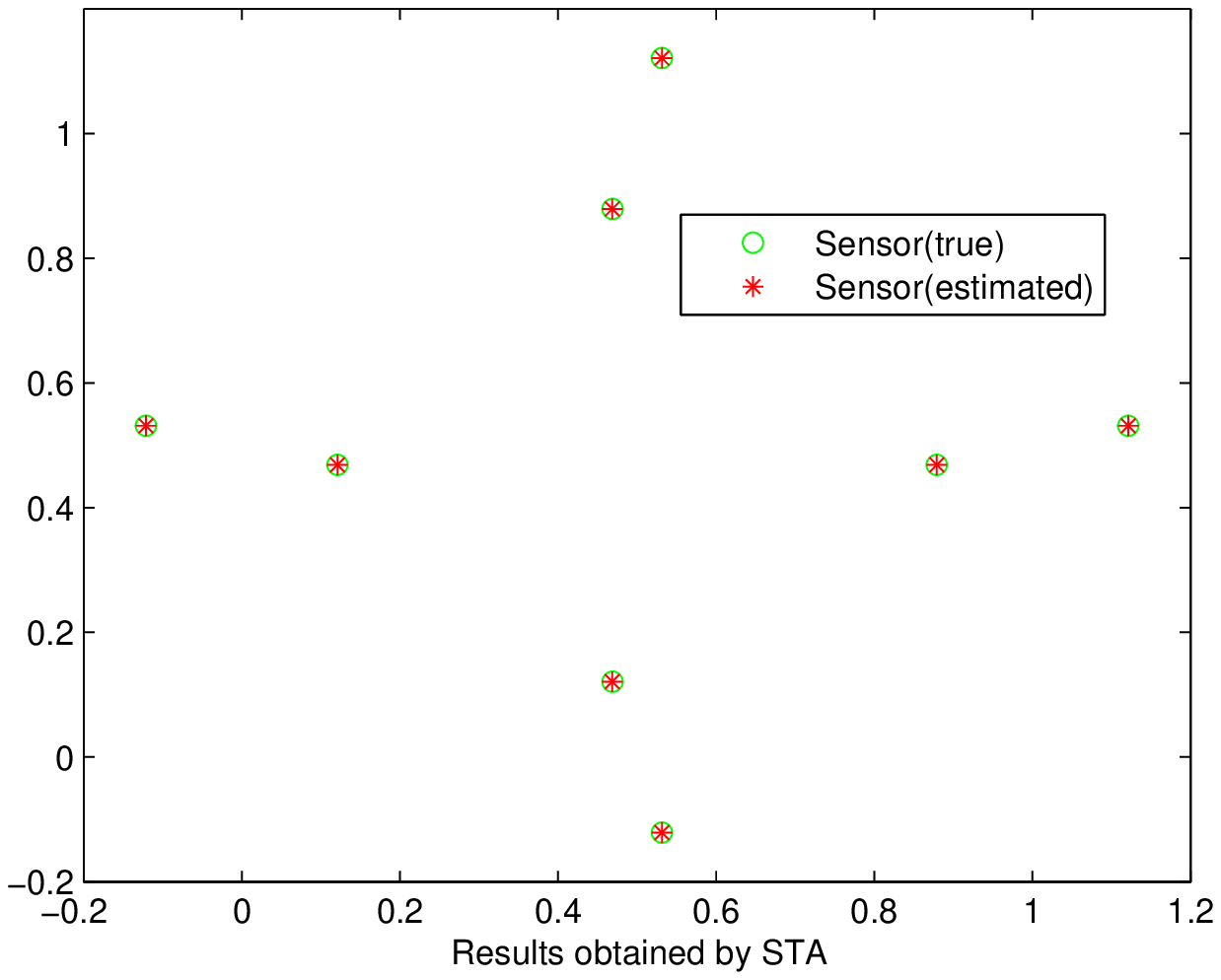} \\
  \includegraphics[width=8cm]{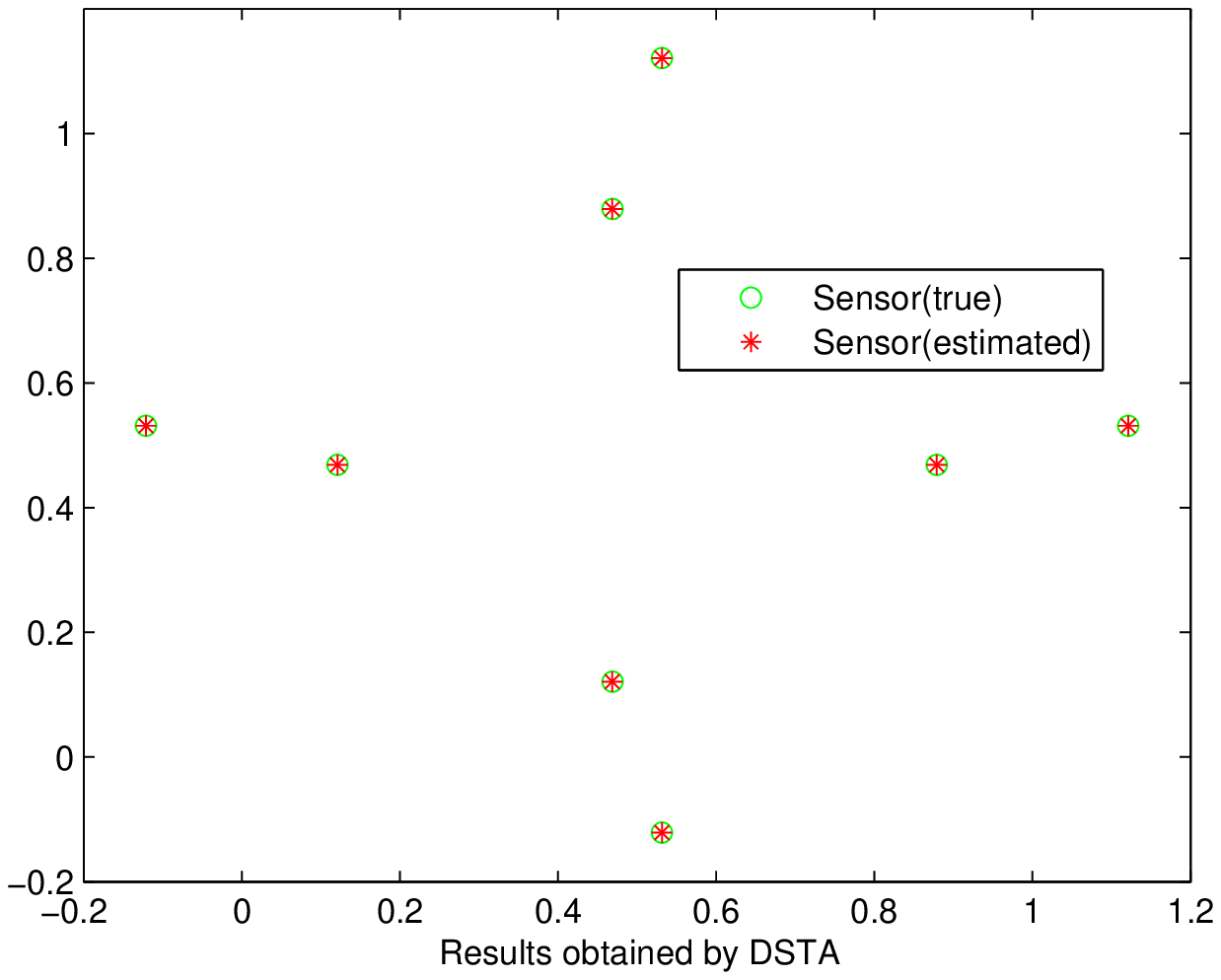}
  \includegraphics[width=8cm]{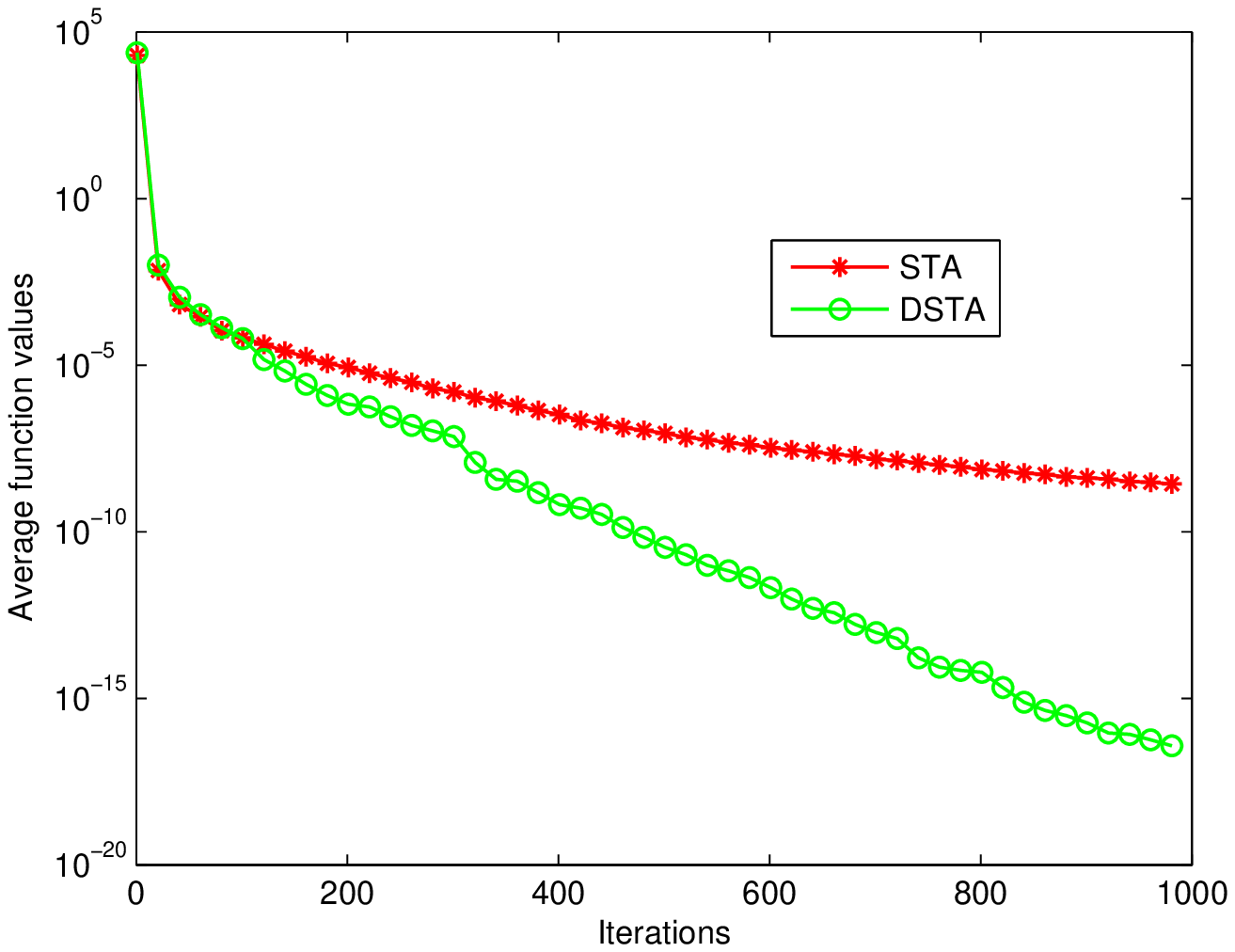}
  \caption{Performance comparison of three different algorithms for the illustrative example}
  \label{fig_illustrative_results}
\end{figure}

For the illustrative example, the parameters setting for the state transition algorithm (STA) and dynamic state transition algorithm (DSTA) is given in Table \ref{parameters_setting}. The search enforcement $SE = 30$,
the transformation factors are all decreasing exponentially from 1 to $1e$-8 with the base $fc = 2$, and the maximum number of iterations (Maxiter) is 1000.
The programmes that implement the STA and DSTA  are run for 20 trails in MATLAB R2010b on Intel(R) Core(TM) i3-2310M CPU @2.10GHz under Window 7 environment.
To evaluate the performance of the proposed approach, we introduce a sparse version of full semi-definite programming (SFSDP), which is a Matlab package for solving sensor network localization problems, for performance comparison.
Numerical results are given in Table \ref{wsn_8s4a}. We observe that both STA and DSTA can find the true global solution, while DSTA is better than STA due to its enhanced exploitation ability, which can be verified by the statistics (best, mean, etc) and the iterative curves of the average function values in Fig. \ref{fig_illustrative_results}.
On the contrary, the SFSDP based on gradient cannot find the true global solution.

Next, two larger problem size of the sensor network localization problems are investigated. The first problem contains 50 sensors and 4 anchors. The second contains 240 sensors and 4 anchors. A radio range of 0.3 is used in the study. The locations of the sensors and the anchors are randomly created by the benchmark generator embedded in the Matlab package. In addition, some noise is added to both problems with the noisy factor of 0.001.

\begin{figure}[!htbp]
  \centering
  \includegraphics[width=8cm]{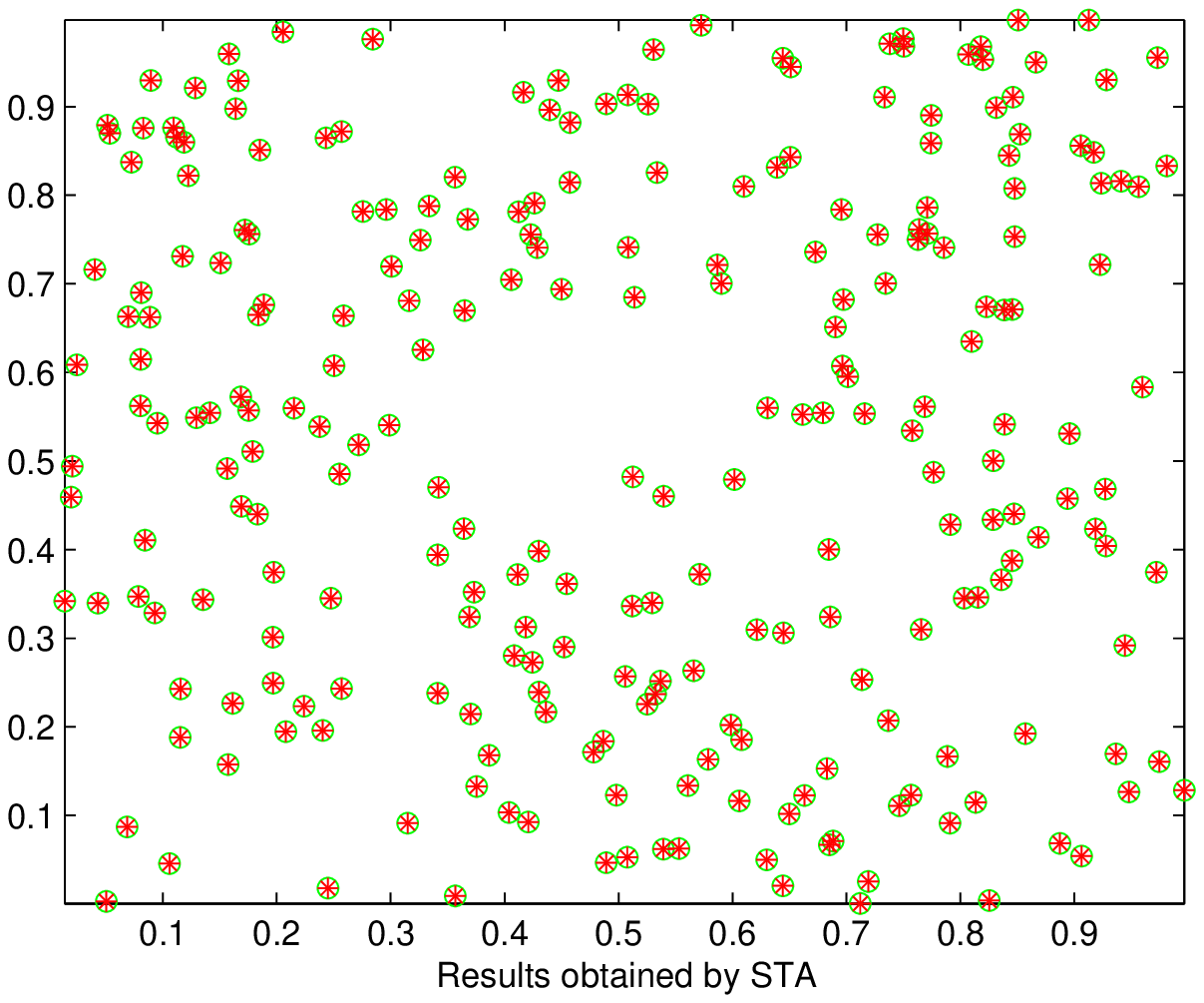}
  \includegraphics[width=8cm]{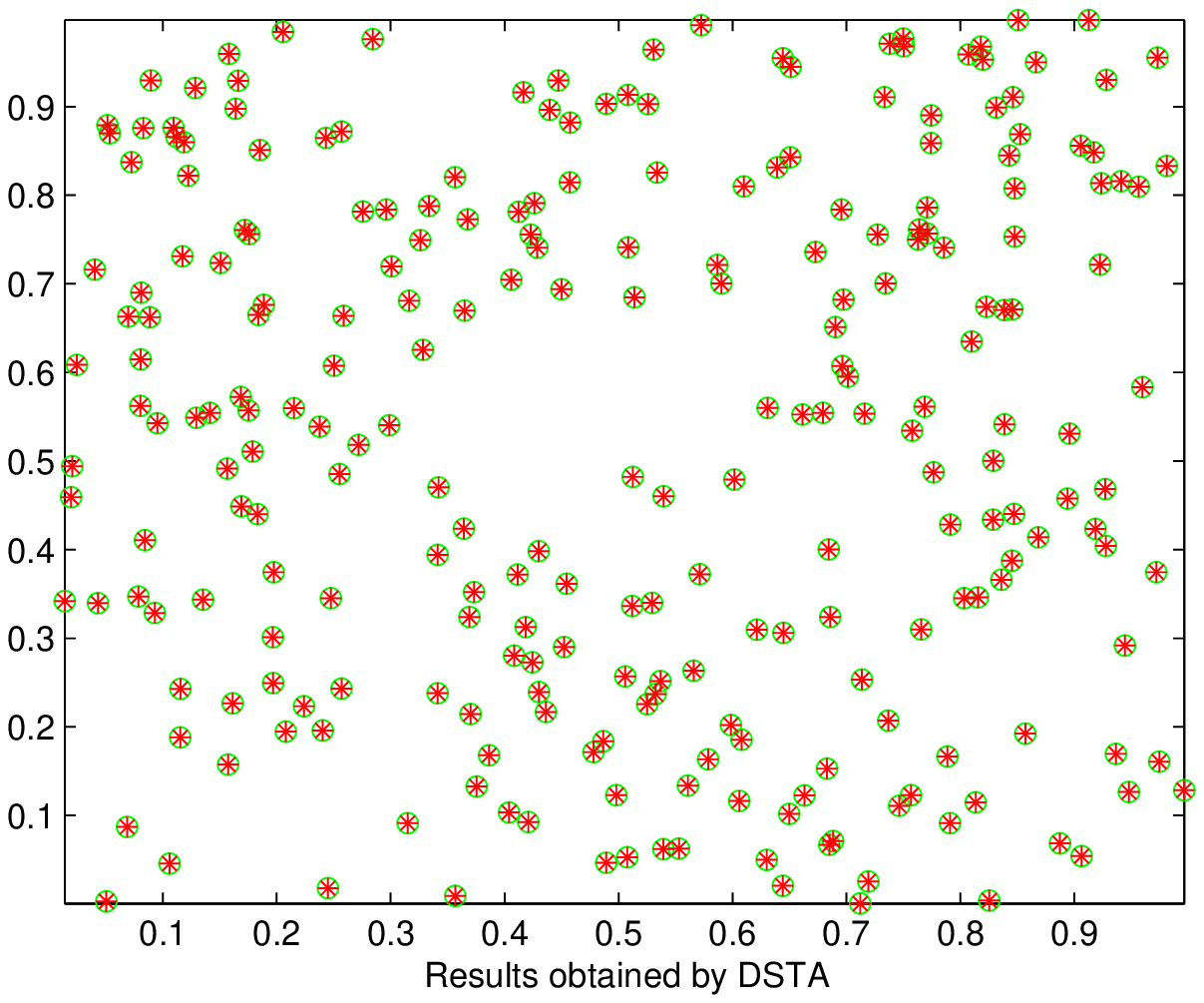}
  \caption{Best results obtained by STA and dynamic STA for the problem with 250 sensors and 4 anchors}
  \label{fig_largesize}
\end{figure}

The experimental results are given in Table \ref{wsn_larger} and Fig. \ref{fig_largesize}. Both STA and dynamic STA can find the global solutions for the two problems, while the performance the
dynamic STA is much better due to its more stable results as indicated by the std.dev.
\begin{table}[!htbp]
\centering
\renewcommand{\arraystretch}{1.3}
\caption{Parameters setting for STA and DSTA}
\scriptsize
\label{parameters_setting}
\begin{tabular}{cccccc}
\hline
Methods  & $\alpha$($\beta$,$\gamma$,$\delta$) & $fc$ & $SE$ & $(p_1,p_2)$ & Maxiter\\
\hline
STA      & $1 \rightarrow 1e$-8   & 2 & 30 &  $-$ & 1000 \\
DSTA      & $1 \rightarrow 1e$-8   & 2 & 30 &  $(0.9,0.3)$ & 1000\\
\hline
\end{tabular}
\end{table}

\begin{table}[!htbp]
\centering
\renewcommand{\arraystretch}{1.3}
\caption{Numerical results for 8 sensors and 4 anchors}
\scriptsize
\label{wsn_8s4a}
\begin{tabular}{{p{3.2cm}|p{3.2cm}|p{3.2cm}|p{3.2cm}}}
\hline
True solutions  & Solutions by STA & Solutions by DSTA & Solutions by SFSDP\\
\hline
    $\bm x_1^{*} = (0.4688,0.1210)$ &  $\bar{\bm x}_1 = (0.4688, 0.1210)$ & $\bar{\bm x}_1 = (0.4688, 0.1210)$ & $\bar{\bm x}_1 = (0.4687, 0.0005)$ \\
    $\bm x_2^{*} = (0.5313,-0.1210)$&  $\bar{\bm x}_2 = (0.5313, -0.1211)$& $\bar{\bm x}_2 = (0.5312, -0.1211)$ & $\bar{\bm x}_2 = (0.5312,  0.0005)$\\
    $\bm x_3^{*} = (0.8790,0.4688)$ &  $\bar{\bm x}_3 = (0.8790, 0.4687)$ & $\bm x_3^{*} = (0.8790,0.4688)$
    & $\bm x_3^{*} = (0.9907,0.4687)$\\
    $\bm x_4^{*} = (1.1210,0.5313)$ &  $\bar{\bm x}_4 = (1.1211, 0.5312)$ & $\bm x_4^{*} = (1.1210,0.5313)$
    &$\bm x_4^{*} = (0.9908,0.5312)$ \\
    $\bm x_5^{*} = (0.5313,1.1210)$ &  $\bar{\bm x}_5 = (0.5312, 1.1211)$ & $\bar{\bm x}_5 = (0.5312, 1.1210)$ & $\bar{\bm x}_5 = (0.5312, 0.9908)$ \\
    $\bm x_6^{*} = (0.4688,0.8790)$ &  $\bar{\bm x}_6 = (0.4688, 0.8790)$ & $\bar{\bm x}_6 = (0.4688, 0.8790)$ & $\bar{\bm x}_6 = (0.4687, 0.9907)$ \\
    $\bm x_7^{*} = (-0.1210,0.5313)$&  $\bar{\bm x}_7 = (-0.1210, 0.5313)$& $\bar{\bm x}_7 = (-0.1210, 0.5312)$ & $\bar{\bm x}_7 = (0.0005, 0.5312)$ \\
    $\bm x_8^{*} = (0.1210,0.4688)$ &  $\bar{\bm x}_8 = (0.1210, 0.4687)$ & $\bm x_8^{*} = (0.1210,0.4688)$ &
    $\bm x_8^{*} = (0.0005,0.4687)$ \\
\hline
\hline
Methods &  Best  &  Mean  & Std. Dev.\\
\hline
STA   & 8.1301e-10 & 2.5723e-9  & 1.4557e-9\\
\hline
DSTA  & 4.2838e-18 & 2.3756e-17 & 2.3502e-17\\
\hline
SFSDP & 0.0171 & 0.0171 & 0\\
\hline
\end{tabular}
\end{table}

\begin{table}[!htbp]
\centering
\renewcommand{\arraystretch}{1.3}
\caption{Numerical results for larger instances}
\scriptsize
\label{wsn_larger}
\begin{tabular}{{p{3.2cm}|p{2.5cm}|p{2.5cm}|p{2.5cm}|p{2.5cm}}}
\hline
Instances  & Methods & Best & Mean & Std. Dev.\\
\hline
50 sensors and 4 anchors & STA  & 2.3129e-5 & 0.0014    & 0.0041\\
                         & DSTA & 2.3129e-5 & 2.3129e-5 & 1.2280e-11\\
\hline
250 sensors and 4 anchors & STA  &  1.8733e-11 & 1.5232e-8  & 5.1607e-8\\
                          & DSTA &  1.6482e-11 & 5.3354e-11 & 3.2959e-11\\
\hline
\end{tabular}
\end{table}

\section{Conclusion}
By using the least squares method, the SNL problem can be reformulated as a non-convex
optimization problem.
In this paper, we have proposed a dynamic state transition algorithm  with refinement for the SNL problem.
A dynamic adjustment strategy called ``risk and restoration in probability" has been incorporated into the basic state transition algorithm for transcending local optimality and improving its global search ability.
Monte Carlo experiments have been designed to obtain a satisfactory combination of the ``risk probability" and
``restoration probability". With the gained parameters setting, the DSTA have been successfully applied to some instances of the SNL problem with a problem size as large as 500.
Numerical results have shown that the proposed DSTA has better performance compared with the basic state transition algorithm in terms of global search ability and solution quality, which has verified the effectiveness and efficiency of the proposed approach.


%

\section*{Acknowledgment}
This work was partially supported by the National Science Foundation for Distinguished Young Scholars of China
(61025015), the Foundation for Innovative Research Groups of the National Natural Science
Foundation of China (61321003), the Australian Research Council (DP140102180, LP140100471) and the 111 Project (B12018).

\ifCLASSOPTIONcaptionsoff
  \newpage
\fi



%
%

\bibliographystyle{ieeetr}
\bibliography{tzy}

%

%






\end{document}